\documentclass[preprint,3p,authoryear,times]{elsarticle}

\journal{}  
\usepackage[fleqn]{amsmath}
\usepackage[]{amsfonts}
\usepackage[]{amsthm}
\usepackage[]{amssymb}
\usepackage{empheq}
\usepackage{lscape} 

\usepackage[mathlines]{lineno}

\usepackage[]{algorithm2e}
\usepackage[english]{babel}
\usepackage[T1]{fontenc}

\usepackage{hyperref} 
\usepackage{graphicx}
\usepackage{caption,subcaption}
\usepackage{tikz}
\usepackage{pstricks}
\usepackage{xstring}
\usepackage{multirow}
\usepackage{hhline}
\usepackage{cases}
\usepackage{xcolor}  

\def\Cc{{\cal C}}

\def\Fc{{\cal F}}
\def\Kc{{\cal K}}

\def\Hc{{\cal H}}
\def\Jc{{\cal J}}
\def\Lc{{\cal L}}

\def\Tc{{\cal T}}

\def\Xc{{\cal X}}

\def\Zc{{\cal Z}}

\def\x{\times}

\def\1{{\bf 1}}

\def \eps{\varepsilon}

\def\argmin{\mathop{\rm argmin}}
\def\argmin_#1{\underset{#1}{\mathrm{argmin\, }}}

\newcommand{\Def}{\stackrel{def}{=}}

\newcommand{\newref}[2]{\hyperref[#2]{#1~\ref*{#2}}} 

\usepackage{longtable}
\usepackage{pdflscape}
\usepackage{rotating,tabularx}
\usepackage[toc,page]{appendix}

\begin{document}
\begin{frontmatter}

\title{Collective Departure Time Allocation in Large-scale Urban Networks: A Flexible Modeling Framework with Trip Length and Desired Arrival Time Distributions}

\author[1]{Mostafa Ameli\corref{cor}}
\cortext[cor]{Corresponding author.} 
\ead{mostafa.ameli@univ-eiffel.fr}
\author[1]{Jean-Patrick Lebacque} 
\author[1]{Negin Alisoltani}
\author[2]{Ludovic Leclercq}

\address[1]{Université Gustave Eiffel, COSYS-GRETTIA, Paris, France}
\address[2]{Université Gustave Eiffel, Université Lyon, ENTPE, LICIT-ECO7, Lyon, France}

\begin{sloppypar} 
\begin{abstract} 
Urban traffic congestion remains a persistent issue for cities worldwide. Recent macroscopic models have adopted a mathematically well-defined relation between network flow and density to characterize traffic states over an urban region. Despite advances in these models, capturing the complex dynamics of urban traffic congestion requires considering the heterogeneous characteristics of trips. Classic macroscopic models, e.g., bottleneck and bathtub models and their extensions, have attempted to account for these characteristics, such as trip-length distribution and desired arrival times. However, they often make assumptions that fall short of reflecting real-world conditions. To address this, generalized bathtub models were recently proposed, introducing a new state variable to capture any distribution of remaining trip lengths. This study builds upon this work to formulate and solve the social optimum, a solution minimizing the sum of all users' generalized (i.e., social and monetary) costs for a departure time choice model. The proposed framework can accommodate any distribution for desired arrival time and trip length, making it more adaptable to the diverse array of trip characteristics in an urban setting. In addition, the existence of the solution is proven, and the proposed solution method calculates the social optimum analytically. The numerical results show that the method is computationally efficient. The proposed methodology is validated on the real test case of Lyon North City, benchmarking with deterministic and stochastic user equilibria.

\end{abstract}


\end{sloppypar} 
\end{frontmatter}


\begin{sloppypar} 

\section{Introduction}
Traffic congestion occurs when the traffic density increases while the traffic flow remains constant or decreases. Macroscopic models aim to rule out urban traffic congestion by holding several assumptions. The common assumption between all macroscopic models is the homogeneity of speed within a single zone based on its traffic density (i.e., accumulation). Seminal work of \cite{geroliminis2008existence} showed that a mathematically well-defined relation between network flow and density could characterize traffic states over an urban region. This concept is very appealing for many applications, including deriving optimal settings for network equilibrium or optimum. Network equilibrium is usually addressed through the concept of User Equilibrium (UE). 
UE in a traffic network, also known as Wardrop's 1$^{st}$ principle, refers to the condition where all commuters select the optimum decision variables (e.g., route or departure time), resulting in no commuters being able to decrease their own travel cost (time) by changing routes \citep{wardrop1952road}. Essentially, all choices have equal and minimal travel times under prevailing traffic conditions. In this context, Stochastic User Equilibrium (SUE), on the other hand, accounts for uncertainty and variability in the cost evaluation (actual vs. perceived cost). As a result, commuters choose where perceived travel costs are minimized, leading to a distribution of travel choices rather than an absolute optimum decision. From the system point of view, authorities aim to minimize the sum of all user social and monetary costs. The solution is called Social or System Optimum (SO), based on the 2$^{nd}$ principle of \cite{wardrop1952road}. The SO solution is the ideal situation for the system, and its calculation and characterization are crucial for improving the transportation system in urban areas \citep{ameli2020improving} as it defines a potential target for authority policies. These principles determine the optimality conditions for decision variables to address the demand. Note that demand characteristics are defined as continuous distributions for all trips or discrete for each trip. The second configuration is known as the trip-based model \citep{mariotte2017macroscopic}. In the context of macroscopic models, a given trip has limited attributes \citep{arnott2022social}.

To calculate SO accurately at the macroscopic level, in order to address the morning commute problem, we need to consider an accurate dynamic model, including the characteristics of trips. The study of peak-hour congestion and congestion models has been ongoing for over 50 years \citep{li2020fifty}. One of the widely used models in macroscopic traffic modeling is Vickrey's bottleneck (point-queue) model, introduced in 1969 by \cite{vickrey1969congestion}, which represented congestion as a fixed-capacity point queue. In classical bottleneck models, A trip is defined by its departure and desired arrival times \citep{li2020fifty}. Therefore, these models involve only the departure time choice dimension for the morning commuters while other travel choice dimensions, such as route, mode and parking choices, and the evening commute, are not considered. 

\cite{luo2020departure} presents a departure time model for estimating the temporal distribution of network-wide traffic congestion during morning rush hours. This model revises the bottleneck model by relaxing the assumption that the last commuter experiences the free-flow travel time. The model is validated using real-world data from Beijing, China. It also proves the effectiveness of the point queue model in estimating travel time, assuming no spillover occurs on road segments. However, this work does not extend the model to other periods like evening rush hours and keeps the assumption of having a homogeneous desired arrival time and a single (averaged) trip length value for all trips, making it not applicable to large-scale heterogeneous urban networks. Additionally, the study does not explicitly calculate user equilibrium or system optimum.

Commuters may also decide on their travel route besides departure time. To address this aspect, Macroscopic fundamental diagram (MFD) or bathtub models have been developed, which take into account trip-length distribution in addition to departure time and desired arrival time distributions. The classic bathtub model, proposed by \cite{vickrey1969congestion}, defines the network as an undifferentiated movement area with a mean speed function that decreases as demand increases \cite{arnott2013bathtub}. MFD models followed the same concept. Further advancements are done in the literature, including considering different desired arrival times \citep{fosgerau2015congestion} for the bathtub model and incorporating trip length distributions through the trip-based MFD \citep{mariotte2017macroscopic, leclercq2017dynamic}. However, most models extending the Classic bathtub model make the assumption that time-independent negative exponential distribution represents the remaining trip distance of all trips traveling in the system. For example, \cite{arnott2016equilibrium} presented a morning rush-hour traffic dynamics model based on the bathtub model, which incorporates hypercongestion situations of heavy congestion where throughput decreases as traffic density increases. Despite the importance of hypercongestion in real-world traffic dynamics, it has been challenging to incorporate it into models due to analytical intractability. The authors developed a simplified model that allows for some degree of analytical tractability, enabling them to study the properties of the model under equilibrium and optimal conditions. The configuration of the mentioned study focused on an isotropic downtown area with identical commuters, utilizing Greenshields' simplified fundamental  diagram and a specific cost function. However, the model's assumptions, particularly the fixed departure times and absence of late arrivals, limit its flexibility in capturing more realistic variations in commuter behavior. This rigid departure time choice model may not sufficiently account for the complexities and variability of commuter behaviors in response to changes in traffic conditions and work schedules. 

\cite{amirgholy2017modeling} formulated the dynamics of congestion in large urban networks using the MFD and examined the morning commute problem. They developed a bathtub model by combining Vickrey's model of dynamic congestion with the MFD to formulate the user equilibrium. The paper presented both exact numerical solutions and analytical approximations of the user equilibrium condition. Moreover, it proposed dynamic tolling and taxing strategies to minimize the generalized cost of the system. However, their approach heavily relied on a well-behaved remaining trip-length function for analytical approximation. 

The assumption regarding the remaining trip distance in classic MFD/bathtub models is not representative of real-world test cases, as shown by different empirical studies, e.g., \cite{liu2012understanding}. This assumption is necessary because the state variable used to capture the dynamics in these models is always accumulation, i.e., the number of users in the network at time $t$, in the classic MFD/bathtub models \citep{laval2022effect}. Several studies in the literature have extended the classic bathtub model to incorporate heterogeneous trip length distributions \citep{lamotte2018morning} and supply profiles \citep{mariotte2017macroscopic, leclercq2017dynamic} with the same state variable. 

\cite{lamotte2018dynamic} introduced the M-model with the total remaining travel distance as a state variable to provide a computational approximation of the trip-based model. They validated their methodology through numerical experiments using real and simulated data. More recently, \cite{jin2020generalized} proposed the generalized bathtub model that extends the classic bathtub model to capture various distributions of the trip length by introducing a new state variable: the number of active trips at time $t$ with remaining distances greater than or equal to threshold $x$, denoted by $K(x, t)$. He formulates the traffic dynamics by four equivalent partial differential equations that track the distribution of the remaining trip lengths. \cite{laval2022effect} investigated the impact of trip-length distribution on the accumulation variance of different macroscopic models and showed that the generalized bathtub model results are valid in both cases of slowly- or rapidly-varying demand. This study aims to formulate and solve the system optimum, also known as the social optimum, for the departure time choice model based on the generalized bathtub model with heterogeneous trip attributes and a generic form of the objective function. In particular, the proposed framework can address any distribution for desired arrival time and trip length.

The literature on departure time choice in the context of macroscopic Dynamic Traffic Assignment (DTA) frameworks is limited \citep{aghamohammadi2020dynamic}. \cite{zhong2021dynamic} conducted a study on dynamic user equilibrium for departure time choice in a trip-based model in an isotropic urban network. While the paper provides a detailed investigation of the dynamic user equilibrium, one of its primary limitations is the assumption of identical travelers, which may not hold in diverse urban settings. To address the demand heterogeneity, basically, numerical methods have been developed to compute departure time distributions and resolve equilibrium conditions \citep{arnott2018solving, lamotte2018congestion} rather than considering the SO conditions.

For the calculation of UE/SUE or SO using MFD/bathtub models, multiple studies in the literature have explored the possibility of relaxing the homogeneity of the trip's characteristics \citep{loder20193d, sirmatel2021modeling, zhong2021dynamic, bao2021leaving, guo2023day}. However, in these studies, there was always an assumption that at least one attribute of the travelers' trips (such as trip length or desired arrival time) is either identical or uniformly distributed. This restricts the optimal departure pattern to specific assumption that may not hold for some real urban cases for the optimal departure time distribution. For instance, \cite{fosgerau2015congestion} presented compelling results showing that under "regular sorting," where shorter trips depart later and arrive earlier compared to longer trips, the problem simplifies significantly for a single MFD, reducing the need for explicit computation of reservoir dynamics. However, \cite{lamotte2018morning} contradicted these findings by demonstrating that a First-in, First-out (FIFO) sorting pattern emerges within user groups with similar scheduling preferences but different trip lengths when there is a single peak in the morning commute. \cite{ameli2022departure} addressed these discrepancies for the UE problem using generalized bathtub models in order to resolve the debate by showing that there is no such property for the equilibrium solution with a fully heterogeneous demand profile.

As mentioned, recent studies have explored generalized bathtub models and applied Mean Field Game theory for deterministic user equilibrium \cite{ameli2022departure} and Stochastic User Equilibrium analysis \citep{lebacque2022stochastic, ameli2023morning}. However, there is a need for further research to address the social optimum for the morning commute problem within these modeling frameworks. Regarding other mentioned models, recently, the SO problem has been well defined and addressed just for the classic bathtub model by \cite{arnott2022social}. \cite{aghamohammadi2020dynamic}, in their review paper, mentioned necessitates of further research on establishing analytical solutions for the SO conditions, resolving discrepancies, and refining capacity constraints to enhance the understanding and modeling of traffic dynamics within the DTA framework. Indeed, these are the ultimate goals of this study. 

\begin{sidewaystable}[]
\centering
\resizebox{\textwidth}{!}{\begin{tabular}{l|cccc|ccc|ccc|ccc|cc|cc}
\multicolumn{1}{c|}{\multirow{4}{*}{Research}}     & \multicolumn{4}{c|}{Macroscopic Model}                                                                                                                                                                                                                              & \multicolumn{3}{c|}{Equilibrium Formulation}                        & \multicolumn{6}{c|}{Demand Profile}                                                                                                                                                                                                                                                                                                                                                                                                                                      & \multicolumn{2}{c|}{Problem configuration}               & \multicolumn{2}{c}{Optimization method}              \\ \cline{2-18}
                           & \multirow{3}{*}{\begin{tabular}[|c|]{@{}c@{}}Point \\ Queue  \end{tabular}} & \multirow{3}{*}{\begin{tabular}[|c|]{@{}c@{}}MFD \\ (NFD)  \end{tabular}} & \multirow{3}{*}{\begin{tabular}[|c|]{@{}c@{}}Classic \\ Bathtub  \end{tabular}} & \multirow{3}{*}{\begin{tabular}[|c|]{@{}c@{}}Generalized \\ Bathtub  \end{tabular}} &  \multirow{3}{*}{\begin{tabular}[|c|]{@{}c@{}} UE \end{tabular}}  &  \multirow{3}{*}{\begin{tabular}[|c|]{@{}c@{}} SUE \end{tabular}} &  \multirow{3}{*}{\begin{tabular}[|c|]{@{}c@{}} SO \end{tabular}}  & \multicolumn{3}{c|}{Trip length}                                                                                                                                                                                                      & \multicolumn{3}{c|}{Desired arrival time}                                                                                                                                                                                         & \multirow{3}{*}{\begin{tabular}[|c|]{@{}c@{}} Continuous \end{tabular}} & \multirow{3}{*}{\begin{tabular}[|c|]{@{}c@{}} Discrete \end{tabular}} & \multirow{3}{*}{\begin{tabular}[|c|]{@{}c@{}} Exact \end{tabular}} & \multirow{3}{*}{\begin{tabular}[|c|]{@{}c@{}} Heuristics \end{tabular}} \\ \cline{9-14}
                           &                              &                                                                      &                                                                             &                                                                                &                      &                      &                      & \begin{tabular}[c]{@{}c@{}}Average\\ Value\end{tabular} & \begin{tabular}[c]{@{}c@{}}Exponential \\ Distribution\end{tabular} & \multicolumn{1}{c|}{\begin{tabular}[c]{@{}c@{}}General\\ Distribution\end{tabular}} & \begin{tabular}[c]{@{}c@{}}Single value /\\ Time window\end{tabular} & \multicolumn{1}{c}{\begin{tabular}[c]{@{}c@{}}Uniform\\ Distribution\end{tabular}} & \multicolumn{1}{c|}{\begin{tabular}[c]{@{}c@{}}General\\ Distribution\end{tabular}} &                             &                           &                        &                             \\ \hline

\cite{vickrey1969congestion}, \cite{vickrey2020congestion} & x & &  & & x & &   & x & &  & x & &  & x& & x & \\
\cite{yang1997analysis} & x & &  & & x & &   & x & &  & x& &  & x& & x& \\\cite{lindsey2019equilibrium} & x &  &  &  & x &  &  & x &  &  &  &  &  x& x &  & x &  \\
\cite{li2020fifty} - review paper &x & &  & &  x& x & & x& & & & & x & x&x & x&x \\
\cite{luo2020departure} & x &  &  &  & x &  &  & x &  &  & x &  &  & x &  &  & x \\
\cite{lamotte2021monotonicity} & x &  &  &  & x &  &  & x &  &  & x &  &  & x & x & x &   \\
\cite{wu2021optimization} & x &  &  &  & x &  &  & x &  &  & x &  &  & x &  & x &  \\
\cite{munoz2006system} & x &  &  &  &  &  & x & x &  &  & x &  &  & x &  & x &  \\
\cite{shen2007path} & x &  &  &  &  &  & x & x &  &  & x &  &  & x &  & & x \\
\cite{kuwahara2007theory} & x &  &  &  &  &  & x & x &  &  &  & x &  & x &  & x &  \\
\cite{guo2023day} & x &  &  &  &  &  & x & x &  &  &  &  & x &  & x & x &  \\
\cite{geroliminis2008existence} &   & x &  & & & &  & x & & &  & & & x &  &  &  \\
\cite{geroliminis2009cordon} & &x &  & & x & &   & x & &  & x& &  & &x & & x\\
\cite{lamotte2016morning}  & & x &  & & x & &   & x& x&  & x & x &  & x&x  &x &x \\
\cite{leclercq2017dynamic} &   & x &  & & & & x & & & x & & & & x &  &  & x \\
\cite{mariotte2017macroscopic} &   & x &  & & x& &  & x & & x & & & & x & x &  & x \\
\cite{lamotte2018dynamic} & & x &  & &  & &   & & & x & x & x &  & &x  & &x \\
\cite{lamotte2018morning} & & x &  & & x & &   & & & x & x & x &  & &x  & &x \\
\cite{loder20193d} &  & x &  &  &  & x &  &  & x &  & x &  &  & x &  &  & x \\
\cite{yildirimoglu2021staggered} &  & x &  &  & x& &  & x & & & x & x &  & & x & & x\\
\cite{zhong2021dynamic} &  & x &  &  & x &  &  & x &  &  & x &  &  & x & x &  & x \\
\cite{amirgholy2017modeling} &  & x & x &  & x &  & x & x &  &  &  & x &  &  &  &  & x \\
\cite{vickrey1991congestion}, \cite{vickrey2019types} &   &  & x & & x & &  & & x & & x& & & x &  & x & \\
\cite{liu2012understanding} &   &  & x & & & & & & & x & & & x & x &  &  & x \\
\cite{arnott2013bathtub} & & & x & & x & &   & x & &  & x & &  & x & & x & \\
\cite{fosgerau2015congestion} & & & x & & x & & x  & x & &  & x& &  & x& & &x \\
\cite{arnott2016equilibrium} &  &  & x &  & x &  &x  & x &  &  &  &  & x  &  & x &  & x \\
\cite{arnott2018solving} &  &  & x & &x & &  &  & & x & x & &  &  & x &  & x \\
\cite{bao2021leaving} &  &  & x &  & x &  &  & x &  &  & x &  &  & x &  &  & x \\
\cite{arnott2022social} &   &  & x & & & & x & x &  & & x & & & x & x &  & x \\
\cite{jin2020generalized} &   &  &  & x & & &  & & & x & & x & & x &  &  &  \\
 \cite{laval2022effect} &   &  &  & x & & &  & & & x &  &  & x & x & x &  &  \\
\cite{ameli2022departure} &   &  &  & x & x& &  & & & x & & x& x & x & x &  & x \\
\cite{lebacque2022stochastic} &   &  &  & x & & x &  & & & x & & & x & x &  x &   & x \\
\textbf{This Study} &   &  &  & \textbf{x} & & & \textbf{x} & & & \textbf{x} & & & \textbf{x} & \textbf{x} & \textbf{x} & \textbf{x} & \textbf{x} \\ 
\end{tabular}\label{tab:RG}
}
\end{sidewaystable}

To conclude our literature review, summarize the state-of-the-art (including recent studies), and highlight the contributions of this study, Table 1 illustrates the characteristics of relevant research papers in the literature. The papers are ordered based on their macroscopic model type. As shown in the table, few studies addressed the SO problem, and to the best of our knowledge, no study has formulated SO for the generalized bathtub model. The papers that have no checked symbol for Equilibrium Formulation did not address the equilibrium problem and either propose the dynamic model \citep{vickrey1969congestion, vickrey1991congestion, jin2020generalized} or analyze \citep{geroliminis2008existence, laval2022effect}) or calibrate \citep{lamotte2018dynamic, liu2012understanding}) the Macroscopic models. Besides, This study proposes a methodology to address generic demand profiles. In contrast, all studies on SO problems at least have an assumption on one characteristic of the demand profile, i.e., distributions of trip length or desired arrival time. Furthermore, we derive the functional derivative of total utility with respect to this variation of the departure time profiles, and we show that this functional derivative utility exists and can be explicitly expressed as a functional of the departure profile. We outline an extension of the model to the case where downstream supply restriction is present. Finally, this study presents an SO model in continuous and discrete settings to investigate the mathematical model analytically and apply it to a practical real test case of Lyon North City. 

The remainder of this paper is organized as follows. In the next section, we present the macroscopic model and illustrate how it captures the network dynamics. The SO problem is presented and discussed in Section \ref{Sec:NC}. We also present the solution 
schemes to solve the SO in this section. The studied test case, the numerical experiments, and the results are presented in Section \ref{sec:NE}. This section includes the comparison of the solutions of UE, SUE and SO. Finally, we outline the main conclusions of this paper and mention some future research directions in Section \ref{sec:conclusion}.

\section{Methodology}

The notations are collected in Table \ref{tab:1}. Bathtub models assume that at time $t$, the velocity ($v_t$) is the same for all traveling users. $v_t$ is a function of the network characteristics and the network load, that is to say, the number of travellers in the network at time $t$, $H(t)$. Let us define the characteristic travel distance $z(t)$ as the distance traveled by a virtual user up until time $t$:
\begin{align}
    \label{eq: virtual user}
    z(t) := \int_0^t v_s ds.
\end{align}

\noindent where $v_s=V(H(s))$, and $V$ is assumed bounded from above and below: $0<V_{min}\leq V \leq V_{max}$. Note that $V_{min}$ can be very small but should be $> 0$, thus $z$ is an invertible function, i.e., since $v_t \geq V_{min}>0\ \forall t\in \Tc$,
$z$ is an invertible function. Let $z^{-1}$ denote the inverse function of $z$. Then, we have $z^{-1}\big(z(t)\big) = t$ and $z^{-1}(x)$ represents the time at which the virtual user has reached $x$. Note that the negative exponential distribution of the trip length transforms the generalized bathtub model to the classic bathtub model \citep{jin2020generalized}. Therefore, the results from both models will be identical. It is worth mentioning that the assumption of exponential distribution for the demand profile also transforms other common macroscopic models (e.g., MFD or trip-based MFD models) to the simple accumulation model and results in the same solution (see \citealt{lamotte2018dynamic, laval2022effect} for the details).

Now, let $T(t_d,x)$ denote the travel time of a player departing at time $t_d$ with trip length $x$. Considering Equation \ref{eq: virtual user},  $T(t_d,x)$ can be determined by,
\begin{align}
    \label{eq: trip time}
    T(t_d,x) = z^{-1}\big(x+z(t_d)\big) - t_d.
\end{align}

\begin{table}[!h] 
\begin{center}
\caption{List of notations}
\resizebox{\textwidth}{!}{\begin{tabular}{p{2cm}p{14cm}} \hline
$\Tc$               &   Time horizon.\\
$x$                 &   Vector of trip lengths. \\
$t_a$               &   Vector of desired arrival times.\\ 

$K(x, t)$           &   Number of agents at time $t$ with remaining trip distance greater than $x$.\\
$H(t)$              &   $\coloneqq K(0, t)$. Number of agents at time $t$ in the network.\\
$v(t)$              &   $=V(H(t))$. Common velocity of agents at time $t$.\\
$z(t)$              &   Characteristic travel distance.\\
$T(t_d, x, t)$  &   Travel time of a trip started at $t_d$ with trip length $x$ at time $t$.\\
$m(t_a,x)$          &   Distributions of demand with trip length $x$ and desired arrival time $t_a$.\\
$h(x)$              &  \textcolor{blue}{ Initial accumulation of agents with trip length greater than $x$.} \\
$f(t_a,x,t)$        &   Distributions of departure times $t$ with desired arrival time $t_a$ and trip length $x$. \\



\hline 
\end{tabular}}\label{tab:1}
\end{center}
\end{table}

In departure time choice problems, the travel cost is usually defined based on $\alpha$-$\beta$-$\gamma$ scheduling preferences \citep{fosgerau2015congestion}. That means the cost function is defined as the sum of the travel time and a penalty cost for arriving at $t_d + T(t_d,x)$ instead of the desired arrival time. Specifically, we assume that each player's cost function is given by, 
\begin{align}
\label{eq:cost function}
J(t_d,x, t_a) = \alpha T(t_d,x) + \beta\big(t_a - t_d - T(t_d,x)\big)_+ + \gamma\big(t_d + T(t_d,x) - t_a\big)_+,
\end{align}
where $(y)_+ = \max\{y,0\}$, and $\alpha$ denotes the cost of traveling per unit of time, $\beta$ and $\gamma$ denote, respectively, the cost of earliness and lateness for the traveller arrival. We assume that the travel cost is an increasing function of travel time,  thus $\alpha > \beta$. The dependency of $J$ on $T(x_d,x)$ expresses indirectly the impact of other travellers on a traveller with attributes $t_d$ and $x$.

The cost function defined in Equation \ref{eq:cost function} captures the fact that travellers prefer not to deviate from their desired arrival time (i.e., arrive as close as possible to their desired arrival time) while they do not spend too much time on the traffic.
Note that the dependency of the cost function on
the trip lengths is not emphasized in the notation, while it holds implicitly. 

Let us now complete the description of the bathtub model. The data is given by the distribution of the number of users desiring to arrive at $t_a$ and sharing the same trip length $x$, $m(t_a,x) dt_a dx$ with respect to desired arrival time $t_a\in \Xc_a$ and trip length $x\in\Xc$. The unknown in the SO problem is the distributions of departure times $t\in \Tc$ with desired arrival time $t_a\in\Tc_a$ and trip length $x\in\Xc$. The resulting distribution of traveller departure time pattern is denoted as $f(t_a,x,t)\, dt_a \,dx \,dt$ and could be considered as flow or number of users based on the problem configuration. Thus $f$ satisfies the following convex set of constraints ($\Kc$):
\begin{equation}\label{eq:SystemConstraints}
   (\Kc) \qquad  \left| \;\; \begin{array}{l}
         \int_{\Tc} f(t_a,x,t) dt = m(t_a,x)  \\
         f(t_a,x,t) \geq 0
    \end{array} \right.
\end{equation}
The dynamics of the bathtub system result from the following processes: i) travellers are conserved, ii) travellers travel at speed $v_t=V(H(t))$, iii) travellers exit the system when they have travelled the trip length $x$ (thus yielding the outflow of the system), iv) the travel demand $f(t_a,x,t)$ yields the inflow into the system. The distribution of initial agents with trip length $x$ provides the initial condition of the system. $z(t)$ and $H(t)$ constitute the main dynamic variables. The following set of equations describes the dynamics of the system:
\begin{equation}\label{eq:SystemEquations}
     \left| \;\; \begin{array}{ll}
        z(t) := \int_0^t dt \, V(H(t)) & (\ref{eq:SystemEquations}.1) \\
        H(t) = h(z(t)) + \int_{0}^{t} ds \,\bar{F}(z(t)-z(s),s)  & (\ref{eq:SystemEquations}.2) \\
        \bar{F}(x,t) = \int_x^\infty d\xi \, \int_{\Tc_a} dt_a\, f(t_a,\xi,t)  & (\ref{eq:SystemEquations}.3)
    \end{array} \right.
\end{equation}
Equation (\ref{eq:SystemEquations}.3) defines $\bar{F}(x,t)$ which is the demand at time $t$ of trips with trip-length greater than $x$. Equation (\ref{eq:SystemEquations}.2), which describes the evolution of $H(t)$,  can be understood in the following way. $h(z(t))$ expresses the contribution of the initial travellers present in the system to $H(t)$ whereas the integral $\int_{0}^{t} ds \,\bar{F}(z(t)-z(s),s)$ expresses the contribution of the departure distribution $f$ to $H(t)$, given that the remaining trip length of each travellers diminishes at a rate $V(H(t))$.

Let us assume some regularity conditions on the initial accumulation $h$ and the initial density $k$ ( initial accumulation and density satisfy $h(x)=\int_x^\infty d\xi \, k(\xi)$:  $k$ should be $\in L^\infty (\Xc)$). We assume similar regularity conditions on the data $m$: the travel demand should be $\in L^\infty (\Tc_a \times \Xc \times \Tc) $. In practice the demand should be smooth enough and should not include jumps in user quantities. Further the velocity should be bounded from below, i;e. there exists $V_{min}$ such that : $V(H)>V_{min}>0$, $\forall H$.  Given these regularity assumptions, and assuming that $f\in L^2 (\Tc_a \times \Xc \times \Tc)$, the following results can be established (refer to appendices B,E and G of \cite{ameli2022departure} ):
\begin{itemize}
    \item Equation \ref{eq:SystemEquations} admits a unique solution  with respect to $z$ and $H$ in $\Cc^0 (\Tc)$ space of continuous functions on $\Tc$;
    \item This unique solution depends Lipschitz continuously (and and also weak-continuously) on the initial conditions $k$ and on the demand $f$.
    \item $z^{-1}$, the inverse of $z$, also depends Lipschitz- and weak-continuously on the initial conditions $k$ and on $f$.
\end{itemize}

\section{System Optimum Model} \label{Sec:NC}

\subsection{Formulation and Existence of the solution}\label{Sec:SO-FormulationExistence}

The objective of the SO problem is typically to optimize the total travel cost of travellers. Thus the objective, denoted as $\Jc$, can be viewed as the sum over all travellers costs given by Equation \ref{eq:cost function}, and the $J$s must be calculated using Equation \ref{eq:SystemEquations}. Thus $\Jc$ is given by
\begin{equation} \label{eq:SO-Criterion}
    \begin{array}{l}
         \Jc \Def \int_{\Tc_a\times\Xc\times\Tc} dt_a\,dx\,dt\, f(t_a,x,t) J(t_a,x,t) \\
         \ \qquad \left| \; \begin{array}{ll}
             J(t_a,x,t) = \alpha T(t,x) + \beta\big(t_a - t - T(t,x)\big)_+ + \gamma\big(t + T(t,x) - t_a \big)_+ & (\ref{eq:SO-Criterion}.1) \\
             T(t,x) = z^{-1} \left( x + z(t) \right)  -t & (\ref{eq:SO-Criterion}.2)    \\
             z(.) \mbox{  solution of (\ref{eq:SystemEquations})} &  (\ref{eq:SO-Criterion}.3) 
         \end{array} \right.
    \end{array}
\end{equation}

Actually, $J$ is a function of $t_a,x,t$ through $z$, which itself is a function of $f,h$ through Equation \ref{eq:SystemEquations}. Thus we can also denote $J$ as $J(f,h)$. Referring again to \cite{ameli2022departure} and the appendix therein, it can be shown that $J$ is Lipschitz continuous, and also that it depends Lipschitz- and weak-continuously on the initial conditions and on $f$. These results could be generalized to $f$ chosen in the set of bounded positive measures on $\Tc_a\times\Xc\times\Tc $. In this paper we consider that $f$ belongs to the Hilbert space of square integrable functions $L^2 \left( \Tc_a\times\Xc\times\Tc \right)$, which is sufficiently general and convenient for applications and numerical approximations. Note also that in the definition of $J$ given in Equation \ref{eq:SO-Criterion} we could substitute the block $\beta\big(t_a - t - T(t,x)\big)_+ + \gamma\big(t + T(t,x) - t_a \big)_+ $ with any other suitable convex function $L$. Thus the SO problem can be stated as follows:
\begin{equation}\label{eq:SO-Statement}
    \min_{f\in \Kc} \, \Jc = \int_{\Tc_a\times\Xc\times\Tc} dt_a\,dx\,dt\, f\,  J(f,h)
\end{equation}
with $J(f,h)$ being calculated from Equation \ref{eq:SystemEquations} by equations (\ref{eq:SO-Criterion}.2) and (\ref{eq:SO-Criterion}.3) . 

The convex bounded domain $\Kc$ is closed in $L^2 \left( \Tc_a\times\Xc\times\Tc \right)$, thus also weakly convex, hence weakly compact (refer to subsection 1.3 of \cite{hinze2008optimization}). Given  initial conditions $h$ and data $m$, let us consider a sequence $\{f_n\}_{n\in\mathbb{N}}$ in $(\Kc)$ which minimizes $\Jc$ (since $\Jc$ is bounded from below by 0 such a sequence exists). By weak compacity of $(\Kc)$ in $L^2 \left( \Tc_a\times\Xc\times\Tc \right)$ we can extract a weakly convergent $\{f_{n_k} \}_{k\in\mathbb{N}}$, the limit of which is denoted $f^*$. Since $\Jc$ is weakly continuous with respect to $f$, the limit as $k \rightarrow \infty$ of the sequence $\{\Jc \left(f_{n_k} \right)\}_{k\in\mathbb{N}}$ equals $\Jc(f^*)$ and is the minimum of $\Jc$. Thus Equation \ref{eq:SO-Statement} admits a solution in the functional space $L^2 \left( \Tc_a\times\Xc\times\Tc \right)$. This solution is not necessarily unique. Actually, by a similar argument, we could show the existence of a solution to Equation \ref{eq:SO-Statement} in the set of bounded positive measures. A solution in the measure space should have a better criterion value but exhibit less regularity than a solution in the $L^2$-space. A comment: the existence could also be proven by Weierstrass type arguments, refer for instance to theorem 2.43 in \cite{guide2006infinite}.

Finally note that the optimization should cover a constant period. It means the time horizon, $\Tc$, includes the desired arrival time, and all commuters finish their trip in this period. Otherwise, the boundary should be added to the model. To generalize the formulation, we can include the terminal cost in $J$.
 
\subsection{Gradient of $\Jc$}\label{sec:Gradient}

The main idea for calculating the system optimum Equation \ref{eq:SO-Statement} is the following.
First it can be shown that in the space $L^2 \left( \Tc_a\times\Xc\times\Tc \right)$, $\Jc$ admits a gradient. Recall that $L^2$ is our functional setting for the SO problem. Second, the projector on the convex set $\Kc$ is well-defined and can be numerically calculated in a very efficient way. These two facts pave the way for finding numerical solutions of Equation \ref{eq:SO-Statement} based on projected gradient concepts. Several discretization methods are available, based either on a particle discretization or a cell discretization of Equation \ref{eq:SystemEquations}. We choose a cell discretization for the numerical approximation of Equation \ref{eq:SO-Statement}.

To calculate the SO solution, let us calculate the gradient criterion with respect to the density of the distribution of departure times $f$. We apply a small variation $\delta f$ to $f$ and calculate the corresponding variation of the total travel cost $\delta\Jc$. The Lipschitz continuity of the solutions $H,z$ and of $z^{-1}$ of Equation \ref{eq:SystemEquations} with respect to $f$ \citep{ameli2022departure} shows that if $\delta f$ is small in the $L^2$ sense, then $\delta z, \delta H, \delta z^{-1}$ are small in the $L^\infty$ and $\Cc^0$ sense. The variation $\delta\Jc$ must be expressed as an integral with respect to $\delta f$, thus yielding the gradient $\triangledown\Jc$: 
\begin{equation}\label{eq:SO-grad}
    \delta\Jc = \int_{\Tc_a\times\Xc\times\Tc} dt_a\,dx\,dt\,\, \triangledown\Jc . \delta f
\end{equation}

\noindent To calculate the gradient analytically, we start from the definition of $\Jc$, and we calculate the variation $\delta\Jc$:
\begin{equation}\label{eq:SO-VarJc}
    \delta\Jc = \int J . \delta f\, dt + \int \delta J . f\, dt
\end{equation}

\noindent where the marginal cost is calculated by $\int \delta J . f\, dt$. Note that the definition of $\int \delta J . f\, dt$ is equivalent to calculating congestion duration in the case of the bottleneck model (defined in \cite{vickrey1991congestion, vickrey2020congestion}). In Equation \ref{eq:SO-VarJc}, we need to compute $\delta J$, which can be calculated as follows:

\begin{equation}\label{eq:SO-VarJ}
    \delta J = \delta TA \left[ 1 + \frac{\partial }{\partial TA} \left(\beta\big(t_a - TA\big)_+ + \gamma\big(TA - t_a \big)_+\right) \right]
\end{equation}

\noindent where $TA$ denotes the arrival time distribution, $TA (x,t) = t + T(t,x)$. Now, we need to calculate $\delta TA$. It can be derived from the definition of $TA$, $TA (x,t) = z^{-1} \left( x+z(t)  \right)$. Substracting  $z \left( TA (x,t) \right) =  x+z(t) $ from $ ( z + \delta z ) \left( TA (x,t) + \delta TA (x,t) \right) =  x+z(t) + \delta z(t) $ it follows
\begin{equation*}
    \delta z \left( TA (x,t) + \delta TA (x,t) \right) + z \left( TA (x,t) + \delta TA (x,t) \right) - z \left( TA (x,t)  \right) = \delta z(t) 
\end{equation*}
At the first order approximation 
\begin{equation*}
    \begin{array}{l}
        \delta z \left( TA (x,t) + \delta TA (x,t) \right) = \delta z \left( TA (x,t)  \right)   \\
        z \left( TA (x,t) + \delta TA (x,t) \right) - z \left( TA (x,t)  \right) = V(H(TA(x,t))) . \delta TA (x,t) 
    \end{array}
\end{equation*}
(recall that $z'(t) = V(H(t))$   $ \forall t$). Thus:
\begin{equation}\label{eq:SO-VarTA}
    \delta TA(x,t) = \frac{\delta z(t) - \delta z(TA(x,t))}{V(H(TA(x,t)))}
\end{equation}

\noindent Based on the fact that $z$ is the integral of the velocity of the system (see \ref{eq:SystemEquations}.1), we can express $\delta z$ as follows:

\begin{equation}\label{eq:SO-VarZ}
    \delta z(t) = \int_0^t ds \, V'(H(s))\delta H(s)
\end{equation}

\noindent The variation $\delta H$ can be obtained by the derivative of Equation (\ref{eq:SystemEquations}.2). It is basically defined as the contribution of the initial condition (i.e., $h$) in addition to the contribution of the variation of $f$ (i.e., $\bar{F})$. The contribution of $h$ follows from $h(x) = \int_x^\infty d\xi k(\xi )$:
\begin{equation*}
    \delta h(z(t)) = - k(z(t)) \delta z(t) 
\end{equation*}
In order to evaluate the contribution of $\delta f$ to $\delta H$, let us note first that:
\begin{equation*}
  \partial_x \bar{F}(x,t) = - \int_{\Tc_a} dt_a\, f(t_a,x,t)
\end{equation*}
Now we must evaluate at first order the difference
\begin{equation*}
    \int_0^t ds\, (\bar{F} + \delta \bar{F} ) ( z(t) + \delta z(t) - z(s) + \delta z(s) ) - \int_0^t ds\, (\bar{F} ) ( z(t) - z(s) )
\end{equation*}
This difference is equal at first order to:
\begin{equation*}
    \int_0^t ds\, \int_{t_a\in\Tc_a} dt_a\, \left[  \int_{z(t) - z(s)}^\infty d\xi\, \delta f(t_a,\xi,s)  + \left( \delta z(t) - \delta z(s) \right) . f(t_a,z(t) - z(s), s ) \right]
\end{equation*}

The  expression found for $\delta H (t)$ is then: 

\begin{equation}\label{eq:SO-VarH}
\begin{aligned}
        \delta H (t) = & - k(z(t)) \delta z(t) - \int_0^t ds \int_{\Tc_a} dt_a \left[ f(t_a,z(t)-z(s),s)\right] [\delta z(t) - \delta z(s)]  \\ &  + \int_0^t ds \int_{z(t) - z(s)}^{\infty} d\zeta \int_{\Tc_a} dt_a \, \delta f(t_a,\zeta,s)
\end{aligned}
\end{equation}


\noindent Equation \ref{eq:SO-VarH} seems complicated, but it expresses $\delta H(t)$ as a linear function of the past values of $\delta z$ and $\delta f$, i.e. $s\leq t$ in (\ref{eq:SO-VarH}). 
Let us summarize the previous results in a more concise way.
\begin{itemize}
    \item Equation \ref{eq:SO-VarJ} and Equation \ref{eq:SO-VarTA} can be expressed as
    \begin{equation}
        \delta J =\Lc . \delta z
    \end{equation}
    with $\Lc$ defined as
    \begin{equation*}
        \Lc . \delta z \, (t) \Def \left[ 1 + \frac{\partial }{\partial TA} \left(\beta\big(t_a - TA\big)_+ + \gamma\big(TA - t_a \big)_+\right) \right] . \frac{\delta z(t) - \delta z(TA(x,t))}{V(H(TA(x,t)))}
    \end{equation*} 
    \item Equation \ref{eq:SO-VarZ} is expressed as
    \begin{equation}
        \delta z =\Zc . \delta H
    \end{equation}
     with $\Zc$ defined as
    \begin{equation*}
        \Zc . \delta H \, (t) \Def \int_0^t ds \, V'(H(s))\delta H(s)
    \end{equation*}  
    \item Equation \ref{eq:SO-VarH} is expressed as
    \begin{equation}
        \delta H =\Hc . \delta z + \Fc . \delta f
    \end{equation}
    with $\Hc$ and $\Fc$ defined as
    \begin{equation*}
        \begin{array}{ll}
             \Hc . \delta z \, (t) \Def & - \left[ k(z(t)) + \int_0^t ds \int_{\Tc_a} dt_a \left[ f(t_a,z(t)-z(s),s)\right] \right] . \delta z(t)  \\
             \  & \ \quad +  \int_0^t ds \int_{\Tc_a} dt_a \, f(t_a,z(t)-z(s),s) . \delta z(s) \\  \Fc . \delta f \, (t) \Def & \int_0^t ds \int_{z(t) - z(s)}^{\infty} d\zeta \int_{\Tc_a} dt_a \, \delta f(t_a,\zeta,s)
        \end{array}
    \end{equation*}  
\end{itemize}
The main finding is that now, we obtain $\delta H$ as a function of $\delta z$ and $\delta f$, in addition to Equation \ref{eq:SO-VarZ}, wherein we have $\delta z$ as a function of $\delta H$. We need to eliminate $\delta H$, which requires the numerically straightforward solution of a triangular linear system. Then, by replacing unknowns, respectively, in equations \ref{eq:SO-VarH}, \ref{eq:SO-VarZ}, \ref{eq:SO-VarTA}, and \ref{eq:SO-VarJ}, we can express $\delta J$ in terms of $\delta f$ and calculate the $\delta \Jc$ in Equation \ref{eq:SO-VarJc}. 
In operator terms, the following relationships result:
\begin{equation}\label{eq:SOVarComplete}
    \begin{array}{rcl}
       \delta H  & = &  \left(  I - \Hc \Zc  \right )^{-1} \Fc \delta f  \\
       \delta z  & = &  \Zc \left(  I - \Hc \Zc  \right )^{-1} \Fc \delta f  \\  
       \delta J  & = &  \Lc  \Zc \left(  I - \Hc \Zc  \right )^{-1} \Fc \delta f   
    \end{array}
\end{equation}
$I$ denotes the identity. Note that all operators are bounded, as a result of the regularity properties of the solutions of Equation \ref{eq:SystemEquations}. Operators $\Zc$ and $\Hc$ are triangular in the sense that $\Hc . \delta z \, (t)$  and $\Zc . \delta H \, (t)$ only depend on past values of $\delta z (s)$ and $\delta H (s)$, i.e such that $s\leq t$. Finally, $\triangledown \Jc$ is derived from Equation \ref{eq:SO-grad}:
\begin{equation}\label{eq:Gradient}
    \nabla \Jc = J + \Fc' \left(  I - \Zc' \Hc'  \right )^{-1} \Zc' \Lc' . f
\end{equation}
In Equation \ref{eq:Gradient} all operators depend on $f$ (via $z$ and $H$). If $A$ is an operator $A'$ denotes the transpose of $A$. The marginal costs result from Equation \ref{eq:Gradient} and are expressed as $\Fc' \left(  I - \Zc' \Hc'  \right )^{-1} \Zc' \Lc' . f$. They include congestion costs and arrival time penalties. The complexity of Equation \ref{eq:Gradient} results from the fact that although all travellers have the same velocity they do not have the same trip lengths and desired arrival times. The operator $\Fc' \left(  I - \Zc' \Hc'  \right )^{-1} \Zc' \Lc' $ expresses the impact of any $(t_a,x,t)$ category of travellers (in terms of their departure time density $f(t_a,x,t)$ ) on any other $(t_a',x',t')$ category. In a discretized setting this operator is approximated by a matrix.

\subsection{Calculation of the system optimum}\label{sec:SO-Calculation}
The first order optimality conditions for the system optimum Equation \ref{eq:SO-Statement}  can be expressed as
\begin{equation}\label{eq:SO-1stOrderOPtimality}
    f = P_{\Kc} \left[ f + \vartheta \nabla\Jc \left( f , h \right)  \right] \quad \forall \vartheta > 0
\end{equation}
These 1st order optimality conditions are necessary but not sufficient, since no properties of $\Jc$ guarantying the sufficiency of the 1st order conditions Equation \ref{eq:SO-1stOrderOPtimality} can be demonstrated (such as convexity of $\Jc$). Therefore the occurrence of local optima cannot be excluded.

The calculation of the gradient and the projector on $\Kc$ allows us to use any projected gradient-like algorithm. Typically an iterative projected gradient algorithm for solving Equation \ref{eq:SO-Statement} can be formulated as
\begin{equation}\label{eq:PrjectedGradient}
    f^{\tau +1} = P_{\Kc} \left[ f^{\tau} + \vartheta^\tau \triangledown\Jc \left( f^\tau , h \right)  \right]
\end{equation}
where $\tau$ denotes the iteration index, $P_{\Kc}$ denotes the projector on $(\Kc )$, and $\vartheta^\tau$ denotes a coefficient to be adjusted in order to guarantee the decrease of $\Jc$. The divergent series rule provides a simple choice of $\vartheta^\tau$ ($\lim_{\vartheta^\tau \rightarrow \infty} = 0 $, $\sum_{\tau} \vartheta^\tau = +\infty$. This rule yields satisfactory results in numerical tests. The discretization of the calculation of $\triangledown\Jc \left( f , h \right) $ can be carried out based on Equation \ref{eq:SOVarComplete}, Equation \ref{eq:Gradient} and Equation \ref{eq:PrjectedGradient}. 

An alternative method of calculation of the SO would be to use the marginal costs, that is to calculate the user optimum with costs given by Equation \ref{eq:Gradient}, i.e the costs $J + \Fc' \left(  I - \Zc \Hc'  \right )^{-1} \Zc \Lc' . f$. The algorithm Equation \ref{eq:PrjectedGradient} converges towards a local optimum of Equation \ref{eq:SO-Statement}. 

Let us give now the principle of the calculation of the projector $P_{\Kc}$, which is a continuous bounded operator in $L^2 \left( \Tc_a \times \Xc \times \Tc \right)$. Consider $g \in L^2 \left( \Tc_a \times \Xc \times \Tc \right) $, then $f=P_{\Kc} (g)$ is obtained by solving with respect to $\varphi$ the following optimization problem:
\begin{equation}\label{eq:ProjectorK}
    \begin{array}{l}
         \min_{\varphi\in L^2 \left( \Tc_a \times \Xc \times \Tc \right)} | g - \varphi  |^2_{L^2 \left( \Tc_a \times \Xc \times \Tc \right)}   \\
         \ \quad \left| \;\begin{array}{ll}
             \int_\Tc  dt \varphi (t_a,x,t) = m(t_a,x) & \forall t_a\in\Tc_a , x\in\Xc  \\
             \varphi (t_a,x,t) \geq 0 & \forall t_a\in\Tc_a , x\in\Xc , t\in\Tc
         \end{array} \right.
    \end{array}
\end{equation}
Given the linear constraints and the quadratic criterion of Equation \ref{eq:ProjectorK}, the optimality conditions of this program are given by 
\begin{equation}\label{eq:KKT}
    f(t_a,x,t) = P_+ \left[  g(t_a,x,t) + \varsigma (t_a,x) \right]
\end{equation}
where $\varsigma \in L^2 \left( \Tc_a \times \Xc \right)$ is obtained by solving
\begin{equation}\label{eq:KKT-Sigma}
    \int_\Tc  dt \, P_+ \left[  g(t_a,x,t) + \varsigma (t_a,x) \right] = m(t_a,x) \quad \forall t_a\in\Tc_a , x\in\Xc 
\end{equation}
Here $P_+$ denotes the projector on the set of positive numbers, i.e. $P_+ (x) = \max (x,0) $. Since $\varsigma \rightarrow P_+ \left[  g(t_a,x,t) + \varsigma \right]$ is a piecewise linear increasing function of $\varsigma$, Equation \ref{eq:KKT-Sigma} can be solved for all $t_a,x$, yielding the projection $f$ of $g$ on $\Kc$. The method can be easily discretized with the common methods, e.g., cell-wise discretization and particle discretization. 

$\varsigma (t_a,x)$ can be interpreted as the Lagrange coefficient of the constraint $\int_\Tc  dt \varphi (t_a,x,t) = m(t_a,x)$. If there were no positivity constraints, we would obtain the projection of $g$ as $f+\varsigma$, with $\varsigma$ given by $\int_\Tc dx (g + \varsigma ) = m$ i.e. $|\Tc| \varsigma = m - \int_\Tc dx g $. The projector $P_+$ accounts for the positivity constraints that apply to $f$.

Finally let us note that the function $\varsigma (t_a,x) \rightarrow \int_\Tc  dt \, P_+ \left[  g(t_a,x,t) + \varsigma (t_a,x) \right]$ admits a left and right derivative everywhere, and these derivatives are increasing. As a consequence the solution of Equation \ref{eq:KKT-Sigma} can be found numerically by applying a Newton algorithm. 

\subsection{Downstream supply constraint}

In practical instances we may want to apply the model to sub-networks of a large network. Traffic exiting a sub-network is liable to be limited by downstream capacity constraints, say $\sigma(t)$, resulting from downstream congestion. How does such a downstream capacity constraint affect the system Equation \ref{eq:SystemEquations} ? Given the traffic speed $v(t)$ the outflow between $t$ and $t+dt$  is given by $K(0,t)-K(v(t) dt ,t)$. Thus the outflow rate is given by
$$
   v(t) \partial_x K(0,t) = \left( K(0,t)-K(v(t) dt ,t) \right) / dt
$$
The traffic demand of the network can be defined as:
\begin{equation}\label{eq:NetworkDemand}
    \Delta (t)  \Def  - \partial_x K(0,t) .V(H(t))
\end{equation}
If we impose the downstream supply restriction $\sigma (t)$ the outflow rate of the network is the minimum between this supply and the demand $\Delta (t)$. The traffic speed is bounded by this outflow rate which is given by:
$$
 \min \left[ \Delta (t) ,\sigma (t)  \right] 
$$
Thus in presence of a downstream supply constraint the speed of traffic in the network is given by 
\begin{equation}\label{eq:SpeedCongestion}
    v(t) = \min \left[ V(H(t)), \sigma (t) / ( -\partial_x K(0,t))  \right]
\end{equation}
The quantity $-\partial_x K(0,t))$ can easily be calculated. Indeed
$$
   K(x,t) = h_0 (x+z(t)) + \int_0^t \,ds \, \overline{F}(x+z(t)-z(s),s)
$$
(refer to \cite{ameli2022departure}). Since $\overline{F}(x,t) = \int_x^\infty d\xi \, F(\xi,t)$ , it follows
$$
  \partial_x K(x,t) = -k (x+z(t)) - \int_0^t ds \, F(x+z(t)-z(s),s)
$$
and setting $x=0$:
\begin{equation}\label{eq:PracticalDemand}
     -\partial_x K(0,t)) = k (z(t)) + \int_0^t ds \, F(z(t)-z(s),s)
\end{equation}
The speed of traffic is thus given by:
\begin{equation}\label{eq:SpeedSupplyConstraint}
     v(t) = \min \left[ V(H(t)), \sigma (t) / \left( k (z(t)) + \int_0^t ds \, F(z(t)-z(s),s) \right)  \right]
\end{equation}

 Recall that $k$ denotes the initial density with respect to remaining travel distance ($h=-\partial_x h$)  and that $F(x,t) = \int_{t_a \in \Tc_a} dt_a \, f(t_a,x,t)$. Can the system optimum be calculated if there is a downstream supply constraint, i.e. if we apply Equation \ref{eq:SpeedCongestion} and Equation \ref{eq:PracticalDemand} to evaluate the traffic speed $v(t)$? 

The dynamical system (\ref{eq:SystemEquations}) subjected to the downstream capacity constraint $\sigma (t)$ can be expressed as
 \begin{equation}\label{eq:SystemWithDownstreamConstraints}
    \left| \; \begin{array}{ll}
         z(t) := \int_0^t dt \, \min\, \left[ V(H(t)) \, , \, \sigma (s) \large/ \,\left(  k(z(s)) + \int_0^s d\varsigma \, F(z(s)-z(\varsigma),\varsigma)  \right)   \right]  &  (\ref{eq:SystemWithDownstreamConstraints}.1)  \\
         H(t) = h(z(t)) + \int_{0}^{t} ds \,\bar{F}(z(t)-z(s),s)  &  (\ref{eq:SystemWithDownstreamConstraints}.2)
    \end{array}\right.
\end{equation}    
In this formulation the supply constraint is integrated into the dynamical system. 
Thus the SO problem with supply constraint is structurally similar to the SO problem without supply constraints and similar resolution methods should apply.
In order for the ideas of  subsections \ref{Sec:NC}  and \ref{sec:Gradient}  to be applicable to the system (\ref{eq:SystemWithDownstreamConstraints}),  the solution of (\ref{eq:SystemWithDownstreamConstraints}) should exist, be unique and depend continuously on the initial condition $h$ and the demand $f$. Considering (\ref{eq:SystemWithDownstreamConstraints}.1) it appears that more restrictive assumptions must be made on the regularity of $h$ and $f$. Specifically it suffices to assume that $h$ admits a derivative in $L^2 (\Xc)$ and $f$ admits first and second derivatives in $L^2 (\Xc\times\Tc)$ for $h$ and $f$ to be continuous (even Hölder continuous). This property results from the classical Sobolev Embedding Theorem (refer, for instance, to theorem 1.2.26 in \cite{drabek2013methods}). Then, of course, the SO problem must be set in the space $W^{2,2} (\Tc_a \times \Xc\times\Tc)$ of functions which admits first and second derivatives in $L^2 (\Tc_a \times \Xc\times\Tc)$, which is a Hilbert space of the Sobolev type. The calculation of the gradient $\nabla \Jc$ must be adapted accordingly by expressing $\delta \Jc$ with the scalar product of $W^{2,2} (\Tc_a \times \Xc\times\Tc)$.

Further in order for the SO problem to be physically relevant it is also necessary to assume that $\sigma (t)$ is bounded from below, i.e. there exists $\sigma_{min} > 0$ such that $\sigma (t) \geq \sigma_{min}$ $\forall t\in\Tc$. With this condition the velocity $v(t)$ of traffic is bounded from above by $V_{max}$ and also from below by a constant which is $>0$  and which depends on $\sigma_{min} > 0$ and on $f$ (therefore indirectly on the total demand). It follows that both $z$ and $z^{-1}$ are Lipschitz continuous functions of time, and that $T$ and $J$ can be evaluated. Further (\ref{eq:SystemWithDownstreamConstraints}.2) suggests that the dependency of $H$ on $z$ and $f$ is regular (refer to (\ref{eq:SO-VarH}) ). Another way to infer this regularity is to observe that the network inflow does not depend on the velocity and the outflow is either $\sigma (t)$ or an outflow at speed $V(H(t))$ (as in Equation \ref{eq:SystemEquations}). Therefore both network in- and out-flows are regular.
 
With the right choice of a functional space as outlined above, the methods of subsections \ref{Sec:SO-FormulationExistence}, \ref{sec:Gradient} and \ref{sec:SO-Calculation} can be adapted to the SO problem with downstream supply, yielding existence results and numerical methods for the calculation of the SO. These adaptations will be an object of future investigation.

\subsection{Cell-wise discretization of the gradient}

This subsection introduces the principles of the cell-wise discretization of Equation \ref{eq:SO-VarJ}, Equation \ref{eq:SO-VarTA}, Equation \ref{eq:SO-VarZ}, Equation \ref{eq:SO-VarH}, i.e of the operators $\Lc$, $\Zc$, $\Hc$, $\Fc$. Let us consider a discretization of $\Tc, \Xc, \Tc_a$:   $t_n = (n-1)\Delta t$, with $n=1..N$, $x_\ell = (\ell-1)\Delta x$, with $\ell=1..L$ and $ta_{k}$ with $k=1..K$ the set of desired arrival times. A cell $\{ ta_{k} \} \times [x_\ell x_{\ell +1}] \times [t_n,t_{n+1}] $ is denoted cell $(k\ell n)$.  Then $H,z,TA,J, \delta z, \delta_H, \delta TA, \delta J$ are discretized as piecewise linear functions of $t_a, t,x$ characterized by their nodal values
$$
\begin{array}{l}
    z_n = z(t_n), H_n = H(t_n), TA_{\ell n}= TA(x_\ell , t_n), J_{k\ell n}=J(ta_k,x_\ell, t_n), \\
    \delta z_n = \delta z(t_n), \delta H_n = \delta H(t_n), \delta TA_{\ell n}= \delta TA(x_\ell , t_n), \delta J_{k\ell n}=\delta J(ta_k,x_\ell, t_n)
\end{array}
$$
$f$, which is a distribution, is discretized by the following values
$$
   f_{k\ell n} = \int_{t_n}^{t_{n+1}} ds \int_{x_\ell}^{x_{\ell +1}} d\xi \, f(ta_k,\xi,s) 
$$
The unit of $_{k\ell n}$ is the number of passengers. Thus $f$ is discretized as a piecewise constant function the value of which on cell $(k\ell n)$ is $\frac{1}{\Delta t \Delta x}f_{k\ell n}$.
The discretization of the dynamical system Equation \ref{eq:SystemEquations} in this setting has been reported elsewhere (refer for instance to \cite{balzer2023dynamic}). Here we consider only the discretization of the gradient of $\Jc$.

Some operators are easily discretized. For instance Equation \ref{eq:SO-VarZ} yields immediately 
\begin{equation}
    \delta z_n = \frac{\Delta t}{2}  \left[  V'(H_1) \delta H_1 + V'(H_n) \delta H_n  \right] + \sum_{m=2}^{n-1} \Delta t V'(H_m) \delta H_m
\end{equation}
Hence $\Zc_{n,m}$ is equal to 0 if $m>n$, $\Delta t V'(H_m)/2$ if $m=1$ or $m=n$, and $\Delta t V'(H_m)$ if $1<m<n$.

In order to discretize $\delta TA$, it suffices to find $m$ such that $z_m < TA_{\ell n} \leq z_{m+1}$, which given the definition of $TA$ is equivalent to $t_m < x_\ell + t_n \leq t_{m+1}$. Define 
$$\mu_{\ell n}=\left\lceil \frac{x_\ell + t_n}{\Delta x} \right\rceil\, , \quad \alpha_{\ell n}= \mu_{\ell n} - \frac{x_\ell + t_n}{\Delta x} $$
and the following discretization results by Equation \ref{eq:SO-VarTA}:
\begin{equation}
    \delta TA_{\ell n} = \frac{\delta z_n - \alpha_{\ell n} \delta z_{\mu_{\ell n}} + (1- \alpha_{\ell n}) \delta z_{\mu_{\ell n}+1} }{\alpha_{\ell n} V(H_{\mu_{\ell n}}) + (1- \alpha_{\ell n}) V(H_{\mu_{\ell n}+1}) }
\end{equation}
The expression of $\Lc$ follows from Equation \ref{eq:SO-VarJ}. The coefficients of $\Lc$ satisfy $\delta J_{k\ell n} = \sum_m \Lc_{k\ell n , m} \delta z_m $, and are given by 
\begin{equation}
\begin{aligned}
        \Lc_{k\ell n , m} = \left[ 1 + \frac{\partial }{\partial TA} \left(\beta\big(ta_k - TA_{\ell n} \big)_+ + \gamma\big(TA_{\ell n} - ta_k \big)_+\right) \right]  \times  \\
    \left| \; \begin{array}{cl}
       \frac{ 1 }{\alpha_{\ell n} V(H_{\mu_{\ell n}}) + (1- \alpha_{\ell n}) V(H_{\mu_{\ell n}+1}) }  & \mbox{ if } m=n  \\
       - \frac{\alpha_{\ell n} }{\alpha_{\ell n} V(H_{\mu_{\ell n}}) + (1- \alpha_{\ell n}) V(H_{\mu_{\ell n}+1}) } & \mbox{ if } m=\mu_{\ell n}   \\         
        - \frac{1- \alpha_{\ell n} }{\alpha_{\ell n} V(H_{\mu_{\ell n}}) + (1- \alpha_{\ell n}) V(H_{\mu_{\ell n}+1}) } & \mbox{ if } m=\mu_{\ell n} +1  \\
        0 & \mbox{ otherwise }
    \end{array} \right.
\end{aligned}
\end{equation}

\begin{figure}[!h]
\begin{center}
\includegraphics[width=13.5cm]{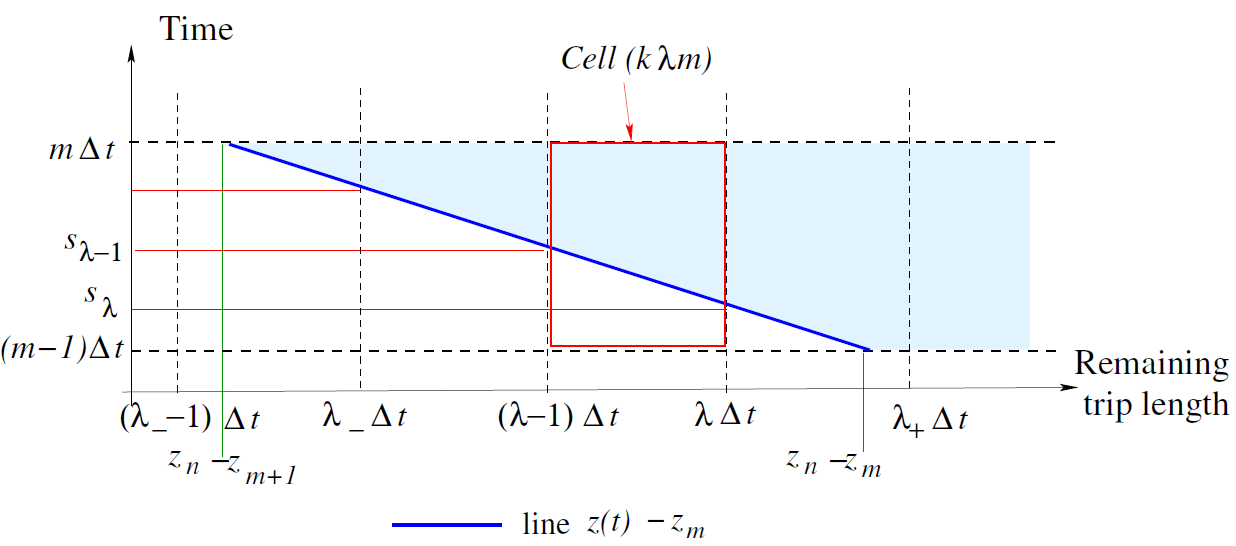}  
\caption{Discretization: calculation of the integrals $E_1$ and $E_2$, contributing to $\Fc$ and $\Hc$.} 
\label{fig:Discretization}
\end{center}
\end{figure}

Other operators are calculated in a similar way. In order to discretize $\Hc$ and $\Fc$ let us consider the following expressions, to be evaluated at $t=t_n$: 
\begin{description}
    \item[i)] $E_1 (t) = \int_0^t ds \int_{z(t) - z(s)}^{\infty} d\zeta \int_{\Tc_a} dt_a \delta f(t_a,\zeta,s)$
    \item[ii)] $E_2 (t) = \int_0^t ds \int_{\Tc_a} dt_a \left[ f(t_a,z(t)-z(s),s)\right] [\delta z(t) - \delta z(s)]  $ 
    \item[iii)] $E_3 (t) = -k(z(t)) \delta z(t) $ 
\end{description}
These expressions contribute to the calculation of $\delta H(t)$ by Equation \ref{eq:SO-VarH}. Let us set $t=t_n$ and discretize $E_{i,n}=E_i(t_n)$ ($i=1,2,3)$. Let us consider $E_1$ first. The integral $\int_0^{t_n} ds $ is replaced by the sum
$\sum_{m=1}^{n-1} \int_{t_m}^{t_{m+1}} ds $. When $s$ varies from $t_m$ to $t_{m+1}$, $z(s)$ varies from $z_m$ to $z_{m+1}$. Define
\begin{equation}
    \lambda_-^{nm} = \left\lceil \frac{z_n - z_{m+1}}{\Delta x} \right\rceil \, , \quad
    \lambda_+^{nm} = \left\lceil \frac{z_n - z_m}{\Delta x} \right\rceil
\end{equation}
The only cells which can contribute to $E_{1,n}$ are the cells $(k\lambda m)$ with $\lambda \geq \lambda_-^{nm}$ as illustrated in figure \ref{fig:Discretization}. For $\lambda > \lambda_+^{nm}$ the contribution is $\delta f_{k\lambda m}$. If $\lambda_-^{nm} \leq \lambda \leq \lambda_+^{nm}$ the contribution is $\alpha^n_{\lambda m} \delta f_{k\lambda m}$, where $\alpha^n_{\lambda m} $ denotes the fraction of the cell $(k\lambda m)$ lying above the curve $z_n - z(s)$ (area shaded in light blue in the figure). This curve is approximated in the plane $(x,t)$ by a straight line joining the points $(z_n - z_{m+1},m\Delta t)$ and $(z_n - z_{m}, (m-1) \Delta t)$. Thus 
\begin{equation}
    E_{1,n} = \sum_{k=1}^K \sum_{m=1}^{n-1} \left[  \sum_{\ell > \lambda_+^{nm}}  \delta f_{k\ell m} + \sum_{\ell = \lambda_-^{nm}}^{\lambda_+^{nm}} \alpha^n_{\ell m} \delta f_{k\ell m}  \right]
\end{equation}
The expression of the operator $\Fc$ can be deduced from
\begin{equation*}
    E_{1,n} = \sum_{k=1}^K \sum_{m=1}^{n-1} \sum_{\ell  \lambda_-^{nm}}^L \Fc_{n,k\ell m} \delta f_{k\ell m} 
\end{equation*}

The evaluation of $E_{2,n}$ follows a similar procedure but requires the evaluation of 
$$
   \sum_{k=1}^K \sum_{m=1}^{n-1}  \int_{t_m}^{t_{m+1}} ds f(ta_k,z_n-z(s),s) . [\delta z_n - \delta z(s)]
$$
The figure \ref{fig:Discretization} illustrates the evaluation of the integral $\int_{t_m}^{t_{m+1}} ds f(ta_k,z_n-z(s),s) . [\delta z_n - \delta z(s)]$. 
The only cells which contribute to this integral are the cells $(k\lambda m)$ with $\lambda = \lambda_-^{nm}, \cdots, \lambda_+^{nm}$. The value of $f$ on each of these cells is $\frac{1}{\Delta t \Delta x}f_{k\lambda n}$. The functions $z_n - z(s)$ and $\delta z_n-\delta z(s)$ are approximated by linear functions of $s$:
\begin{equation}
    \begin{array}{l}
       z_n - z(s) =  z_n - z_m  - \left( -(m-1)+s/\Delta t \right).\left( z_{m+1} - z_m  \right)  \\
        \delta z_n - \delta z(s) =  \delta z_n - \delta z_m  - \left( -(m-1)+s/\Delta t \right).\left( \delta z_{m+1} - \delta z_m  \right) 
    \end{array}
\end{equation}
The integral $\int_{t_m}^{t_{m+1}} ds f(ta_k,z_n-z(s),s) . [\delta z_n - \delta z(s)]$ can now be evaluated in each cell $(k\lambda m)$ with $\lambda = \lambda_-^{nm}, \cdots, \lambda_+^{nm}$ by integrating a linear function of $s$ over the part of the straight line joining the points $(z_n - z_{m+1},m\Delta t)$ and $(z_n - z_{m}, {m-1}\Delta t)$ which lies in each cell $(k\lambda m)$.
\begin{equation}
    \begin{array}{lcl}
       \int_{t_m}^{t_{m+1}} ds f(ta_k,z_n-z(s),s) . [\delta z_n - \delta z(s)] & =  & \sum_{\lambda=\lambda_-^{nm}}^{\lambda_+^{nm}} f_{k\lambda m} \left[  (\delta z_n - \delta z_m) (s_{\lambda +1} - s_\lambda) \right. \\ 
      && \left.\; \; + (\delta z_{m+1} - \delta z_m) \int_{s_\lambda)}^{s_{\lambda +1}} \left( -(m-1)+s/\Delta t \right) ds  \right] \\
       & \Def & \sum_{\lambda=\lambda_-^{nm}}^{\lambda_+^{nm}} \alpha^n_{k\lambda m} \delta z_n  + \beta^n_{k\lambda m}  \delta z_m+ \gamma^n_{k\lambda m} \delta z_{m+1}
    \end{array}
\end{equation}
Thus 
\begin{equation*}
    E_{2,n} = \sum_{k=1}^K \sum_{m=1}^{n-1} \sum_{\ell = \lambda_-^{nm}}^{\lambda_-^{nm}} \alpha^n_{k\ell m} \delta z_n  + \beta^n_{k\ell m}  \delta z_m+ \gamma^n_{k\ell m} \delta z_{m+1}
\end{equation*}

The evaluation of $E_{3,n} = - k(z_n) \delta z_n$ is trivial, it suffices to interpolate $k(z_n)$. Let us define $\mu_-^n = \left\lfloor \frac{z_n }{\Delta x} \right\rfloor$ and $\mu_-^n = \left\lceil \frac{z_n }{\Delta x} \right\rceil$. If $\mu_-^n \Delta x \geq X$ then $E_{3,n} = 0$. If $\mu_-^n \Delta x \leq X < \mu_-^n \Delta x$, $E_{3,n} = \left( \frac{z_n }{\Delta x} - \mu_-^n \right) k_L  \delta z_n$. Finally if $\mu_+^n \Delta x \leq X$ then 
\begin{equation*}
    E_{3,n} = \left( \frac{z_n }{\Delta x} - \mu_-^n \right) k_{\mu_+^n} \delta z_n + \left( - \frac{z_n }{\Delta x} + \mu_+^n \right) k_{\mu_-^n} \delta z_n
\end{equation*}
We can summarize these results by stating that $E_{3,n} = \eps_n \delta z_n$. The coefficients of $\Hc$ result from the identity
\begin{equation*}
    E_{2,n} + E_{3,n} = \sum_{m=1}^n \Hc_{nm} \delta z_m
\end{equation*}
Via $\Hc$ and $\Fc$, $\delta H_n$ is expressed in terms of past values $\delta z_m$ and past values of $\delta f_{k\ell m}$. Only $\delta z_m$ and $\delta f_{k\ell m}$ with values of $m$ such that $m\leq n$ can contribute to $\delta H_n$. Replacing $\Hc.\delta z$ with $\Hc.\Zc \delta H$ in $\delta H =\Hc . \delta z + \Fc \delta f$, $\delta H_n$ can be expressed in terms of a linear combination of past values $\delta H_m$, with $m<1$ and of the values $\delta f_{k\ell m}$ with $m\leq n$. Thus by recursion $\delta H_n$ can be expressed in terms of a linear combination of the $\delta f_{k\ell m}$ only. This means that the discretized version of the operator $\left( I - \Hc\Zc \right)^{-1}\Fc$ can be calculated by matrix products only, without any matrix inversion. The calculation of $\Lc\Zc \left( I - \Hc\Zc \right)^{-1}\Fc$ follows also by matrix products.

%

\subsection{Particle discretization}\label{Sec: particle}

To create a trip-based simulator based on the proposed model, we applied the particle discretization approach, wherein each trip is represented as a particle with multiple attributes, e.g., desired arrival time and trip length. The particle discretization is easily put in correspondence with a micro-simulation. Each particle is endowed with a departure time $TD_p$ (which results from $f$), an arrival time $TA_p$, a remaining trip length $x_p$, and a desired arrival time $t_{a,p}$. The treatment of the dynamics of the system equation \ref{eq:SystemEquations} are now different. The particle $p$ enters the system at time $TD_p$ which is a data,  $x_p$ decreases at a rate $v_n=V(H_n)$, and the particle $p$ exits the system when $x_p=0$ which defines $TA_p$. The number of particles present in the system yields $H$ at any time, thus yield $v$ and $z$. The travel cost for particle $p$ is given by $J_p = \alpha (TA_p - TD_p) + \beta\big(t_{a,p} - TA_p\big)_+ + \gamma\big(TA_p - t_{a,p}\big)_+ $. Finally, the $f$ values follow from $J$ by using the gradient calculation for all particles in a cell. We applied particle discretization to enable our framework to consider trip-based real scenarios, however Cell-wise discretization of the system equation \ref{eq:SystemEquations} yield a faster simulation of the system \citep{nagurney1997projected,ameli2023morning}. 

\section{Numerical Experiments and Results} \label{sec:NE}
In order to prove the concept of the model, the proposed formulation is first applied to a simplified Paris network to assess the performance and effectiveness of the model in a tractable test case. Then, the methodology is applied to a more extensive test case, Lyon North City, in order to evaluate its performance and examine how the optimization procedure for calculating SO affects the congestion level of the network's real state. This application is noteworthy as it is the first instance in the literature that addresses the departure time system optimum on a real large-scale network with a large number of users ant heterogeneous trip profiles (trip lengths and desired arrival time). 

\subsection{Validation of the model: Aggregated scenario of Paris network}
We designed a simplified test case for the center part of Île-de-France. This French region includes Paris City to track the performance of the model in capturing the congestion dynamics and solving the SO problem. The total demand is 1.45 million trips for the morning peak hour. Trips are divided into three classes based on their desired arrival time, and the trip length follows a uniform distribution for two trip length classes. Table \ref{tab:ParisData} presents the characteristics of the demand profile.    

\begin{table}[!h]
\caption{Paris network: Demand profile}
\centering
\begin{tabular}{c|ccc}
Trip length / $t_a$ & 8:00 & 8:30 & 9:00 \\ \hline
{[}0, 18{]}       & 12\% & 20\% & 8\%  \\
{[}18, 42{]}      & 18\% & 30\% & 12\%
\end{tabular}\label{tab:ParisData}
\end{table}

The speed function, $V(.)$, is considered as a piece-wise linear function. The trip-demand has been estimated using the methodological framework provided in \citep{horl2021open}. The parameters of the travel cost are the same for all trips and are chosen based on \cite{ameli2021computational}, which characterizes the travel cost parameters.

\begin{figure}[!h]
\begin{center}
\includegraphics[width=9cm]{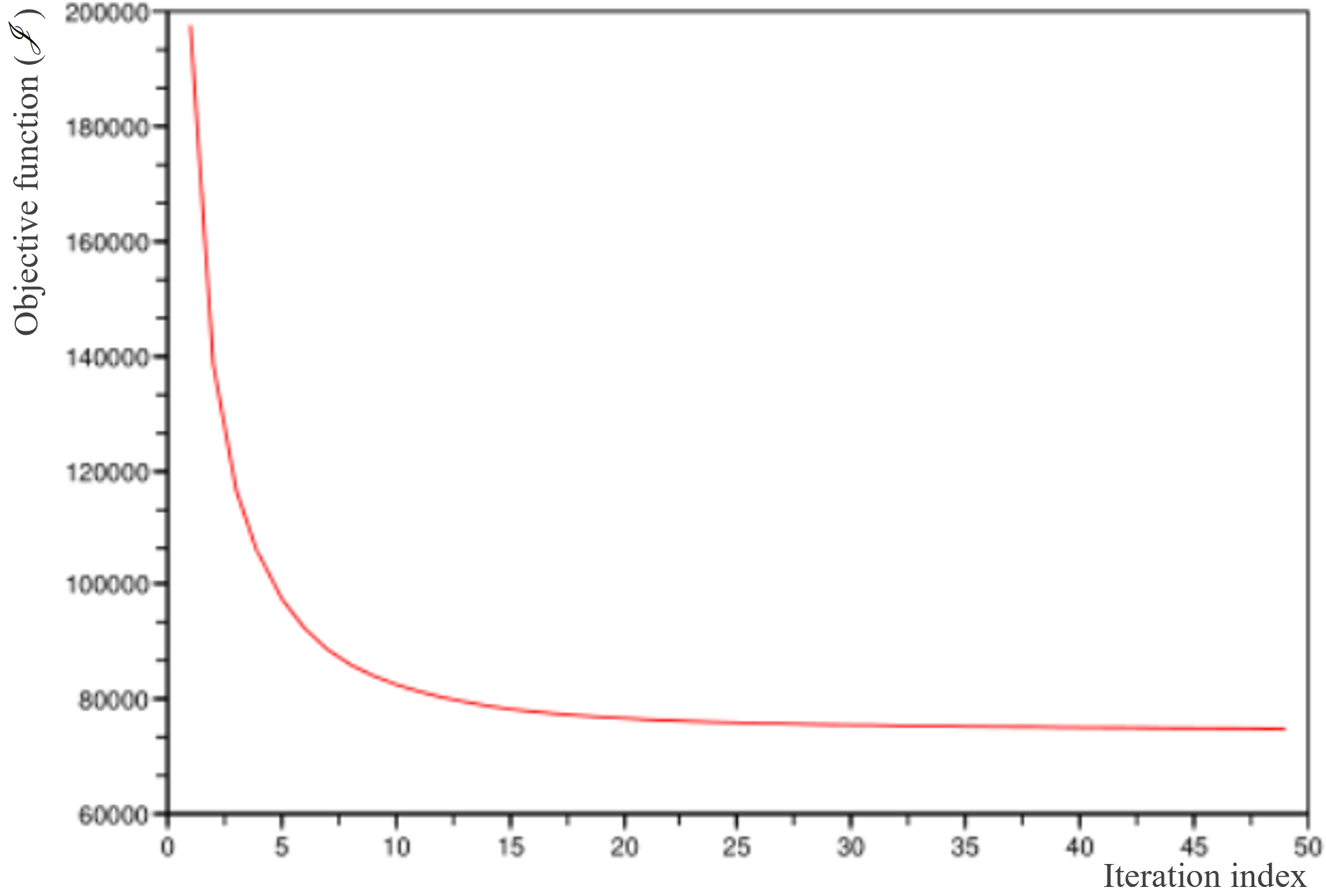}  
\caption{Convergence of the model. $\Jc$ as a function of iteration.} 
\label{fig:Convergence1}
\end{center}
\end{figure}

\begin{figure}[!h]
\begin{center}
\includegraphics[width=8.0cm]{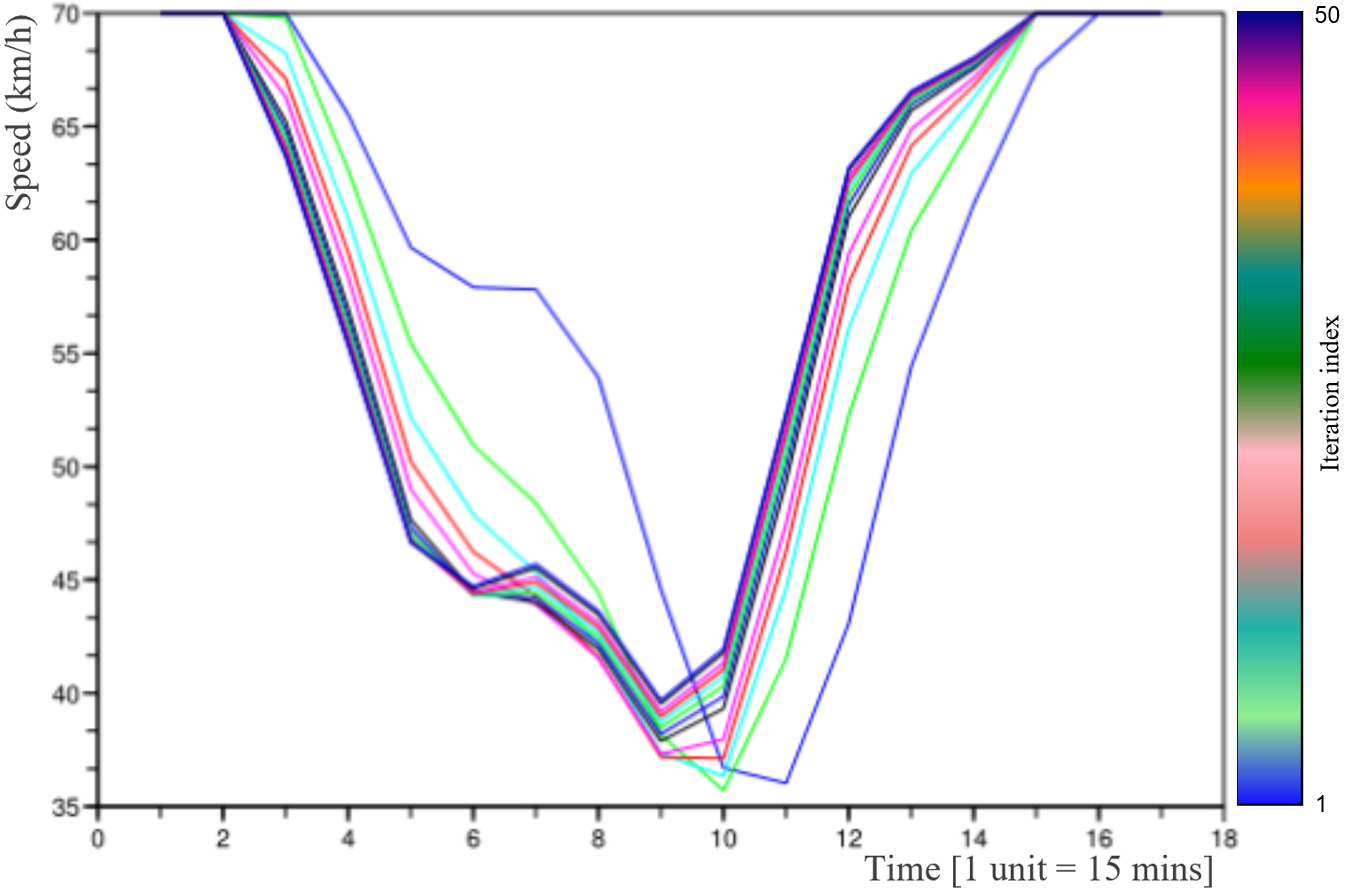}  
\includegraphics[width=8.2cm]{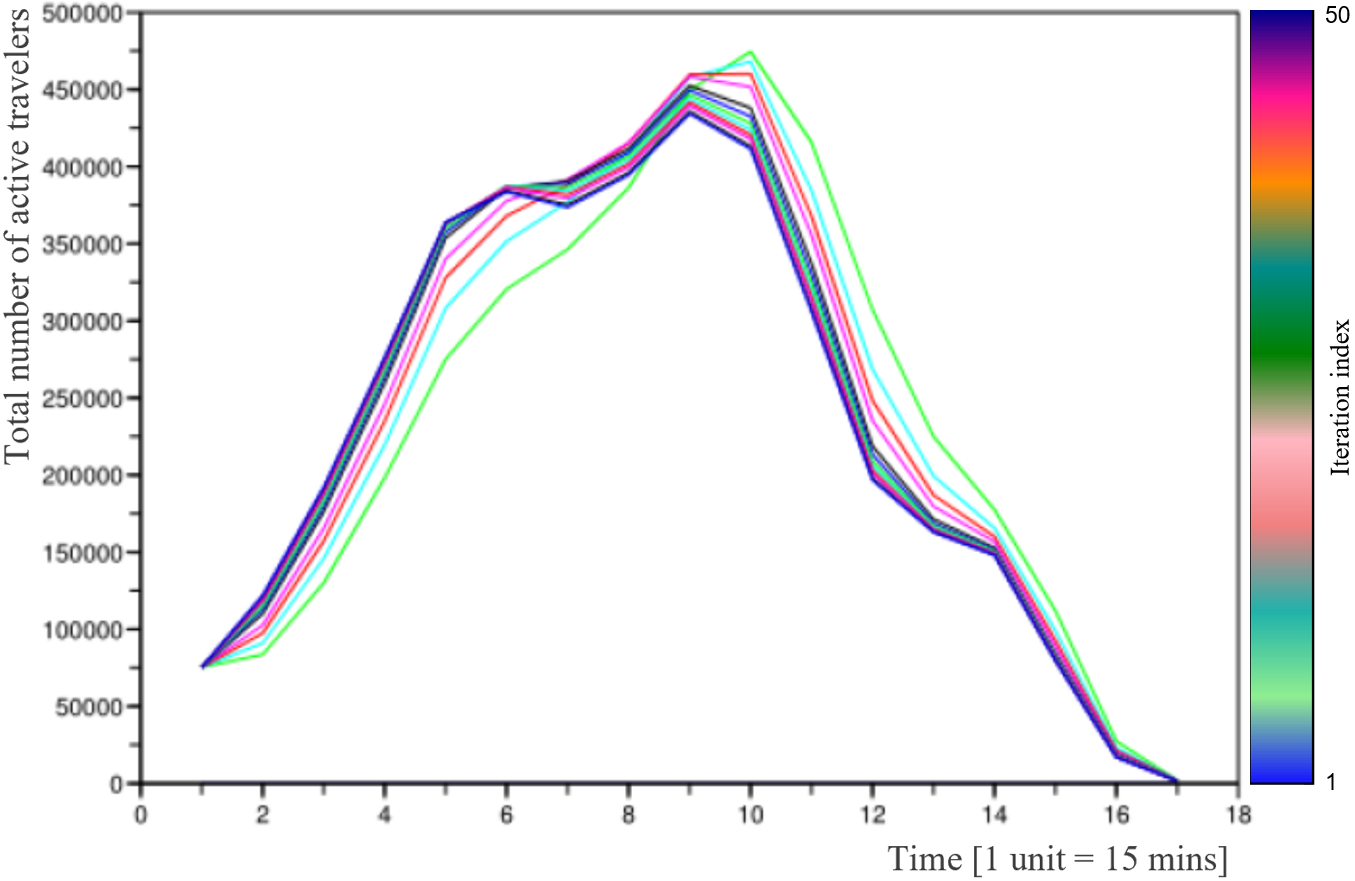}  
\caption{Convergence of the system measures during the optimization. On the left: Speed as a function of time and iteration. On the right: Total number of commuters in the system as a function of time and iteration.} 
\label{fig:itr2}
\end{center}
\end{figure}

 Figure \ref{fig:Convergence1} presents the convergence of the gradient method, presented in Section \ref{sec:Gradient}. It illustrates that the calculation of the gradient leads the algorithm to smoothly converge to the optimum. Close optimal solution obtained with low iteration number, which demonstrates the computational efficiency of the proposed method. Figure \ref{fig:itr2} presents the convergence of the solution method in terms of network criteria. Each curve illustrates the solution of a single iteration. The darkest curve shows the final solution. The left figure presents the network speed as a function of time, and the figure on the right presents the network accumulation as a function of time. In both figures, the final solution has the best value for the targeted criterion (minimum accumulation and maximum network speed). These results prove that the proposed methodology optimizes the network performance following the SO principles.      

\subsection{Large-scale application: Lyon North City}
The proposed SO framework easily adapts to significantly larger instances. In this section, we consider a real-world scenario in the northern region of a French metropolis (Lyon), during the morning peak hours, encompassing 62,450 trips in total. It is worth mentioning that this section illustrates the largest application of our departure time system optimum model to a large-scale network with a realistic demand pattern compared to the literature on departure time choice models.
 \subsubsection{Test case description and demand profile}
\begin{figure}[!t]
\begin{center}
\includegraphics[width=6.3cm]{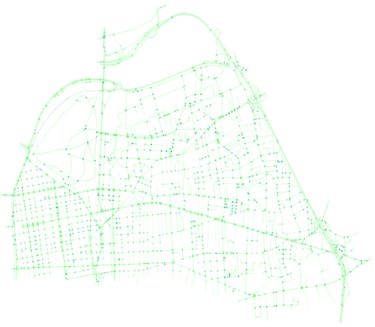}  
\includegraphics[width=8.1cm]{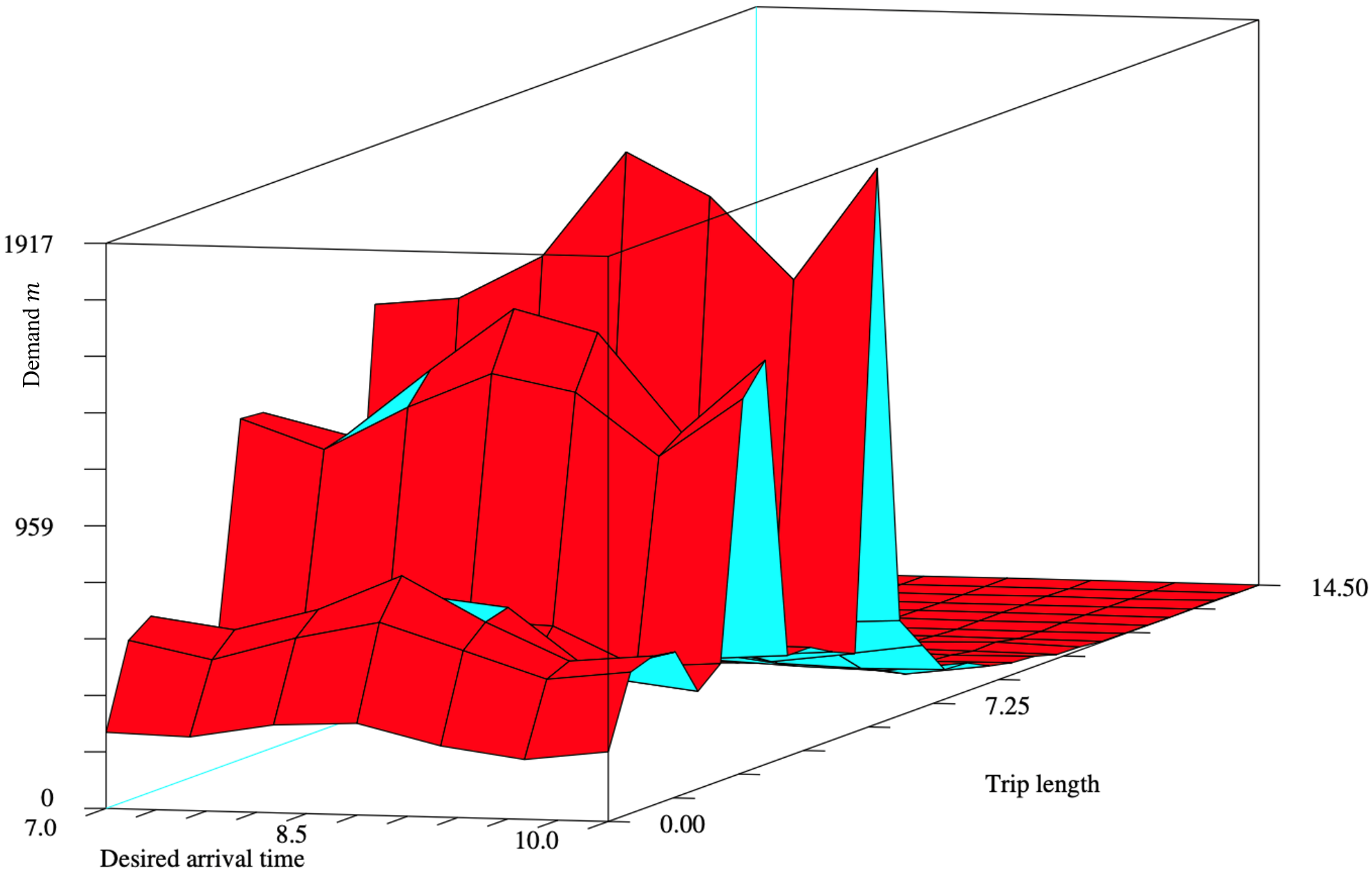}
\caption{The Lyon North data set: On the left: Mapping data \citep{HWN8KE_2021}. On the right: The demand $m$ for the continuous approximation.} 
\label{fig:LyonDemand}
\end{center}
\end{figure}

We implemented and applied our SO model to the Lyon North network, which includes 1,883 nodes and 3,383 links. The network characteristics are presented in \cite{mariotte2020calibration}. The demand profile includes 62,450 trips during the morning peak hours (6:30 AM to 10:40 AM). The data set of the demand profile is published in \cite{HWN8KE_2021}. The network speed function has been calculated in \cite{mariotte2020calibration, alisoltani2022space}. The cost function parameters, i.e., the $\alpha$-$\beta$-$\gamma$ scheduling preferences, are defined based on the study of \cite{lamotte2018morning}: $\alpha = 1$, $\beta = 0.4 + \frac{0.2(k)}{9}$, and $\gamma = 1.5 + \frac{k}{9}$. In order to consider only the heterogeneity of trip length and desired arrival time distributions, $k$ is fixed to $5$ for all trips in this experiment. The resolution is imposed by the original data. The test case has been accurately calibrated to mirror real-world traffic conditions \citep{alisoltani2020multi, alisoltani2021can}. All trips have an origin and destination on the real network and departure times. At the link level of the network (Figure \ref{fig:LyonDemand}), the origin set contains 94 points, and the destination set includes 227 points. This study only retains the original trip lengths, as the generalized bathtub model does not account for local traffic dynamics. Some trips have origins or destinations outside the covered area (51,215 trips) and will not be considered in the optimization of departure time. Note that 11,235 trips are fully interior. The original departure time is disregarded for these, and a desired arrival time is assigned. We categorize them into seven classes with different desired arrival times. The desired arrival time of each user is deduced from the real arrival time of the user, based on real data \citep{ameli2019multi, alisoltani2019optimal}.

\begin{figure}[!h]
\begin{center}
\includegraphics[width=9cm]{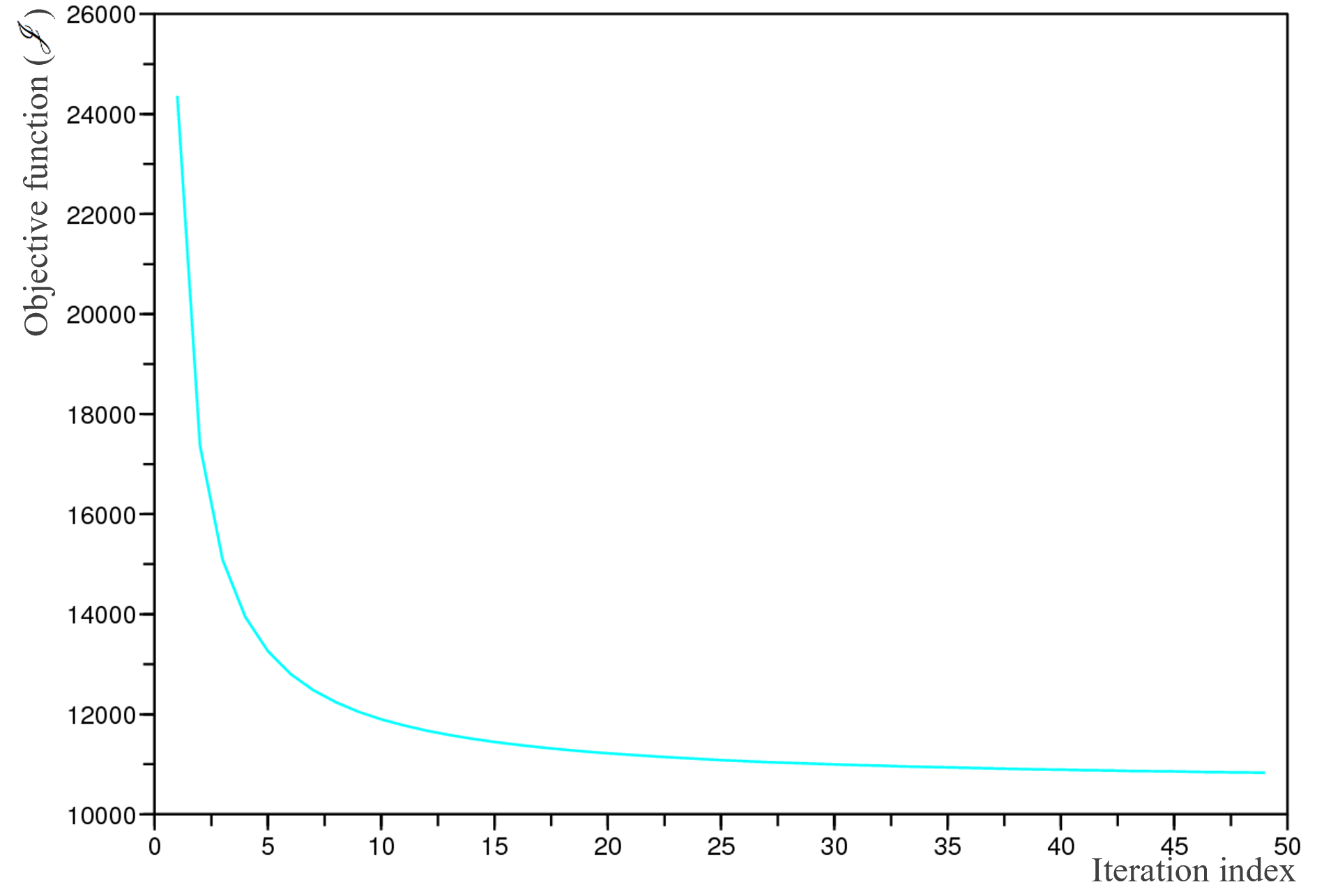}  
\caption{Convergence: SO criterion as a function of iteration.} 
\label{fig:Convergence2}
\end{center}
\end{figure}

Figure \ref{fig:LyonDemand} presents the network graph and the demand profile of this test case. The numerical example was calculated based on a cell discretization (discrete values of $t_a$, cells for $x$ and $t$ values). In the Lyon North case study, 51,215 trips starting or ending outside the study area are excluded from our analysis. This leaves 11,235 completely internal trips. For these, we ignore the original departure times and assign desired arrival times, dividing them into seven categories, starting at 7:30 am, which are separated by half an hour (see Table \ref{tab:Demand}). For this example, convergence is achieved after 25 iterations and finds the optimal solution (see Figure \ref{fig:Convergence2}). Note that the convergence depends on network characteristics' level of congestion.

\begin{figure}[!h]
\begin{center}
\includegraphics[width=8cm]{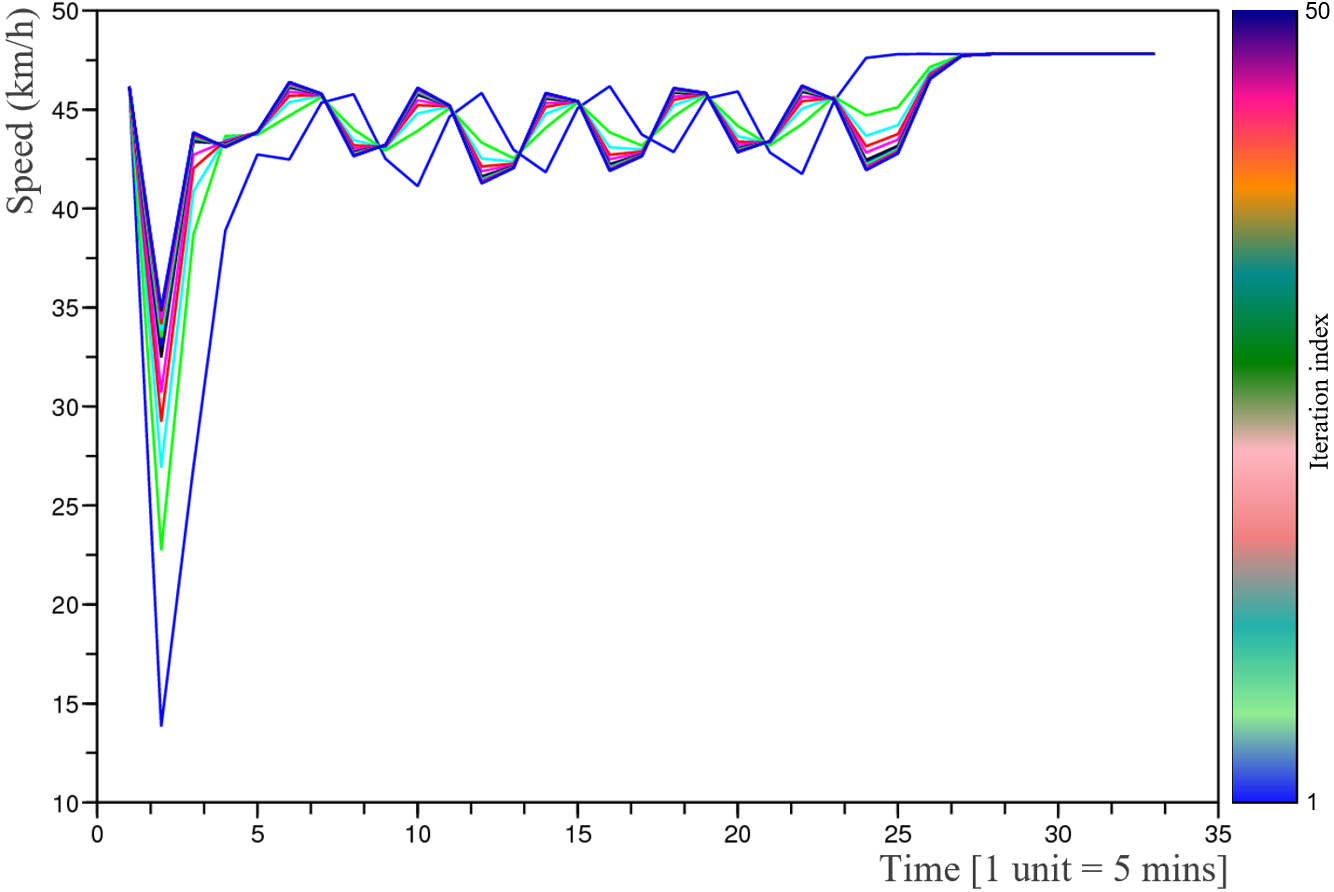} 
\includegraphics[width=8cm]{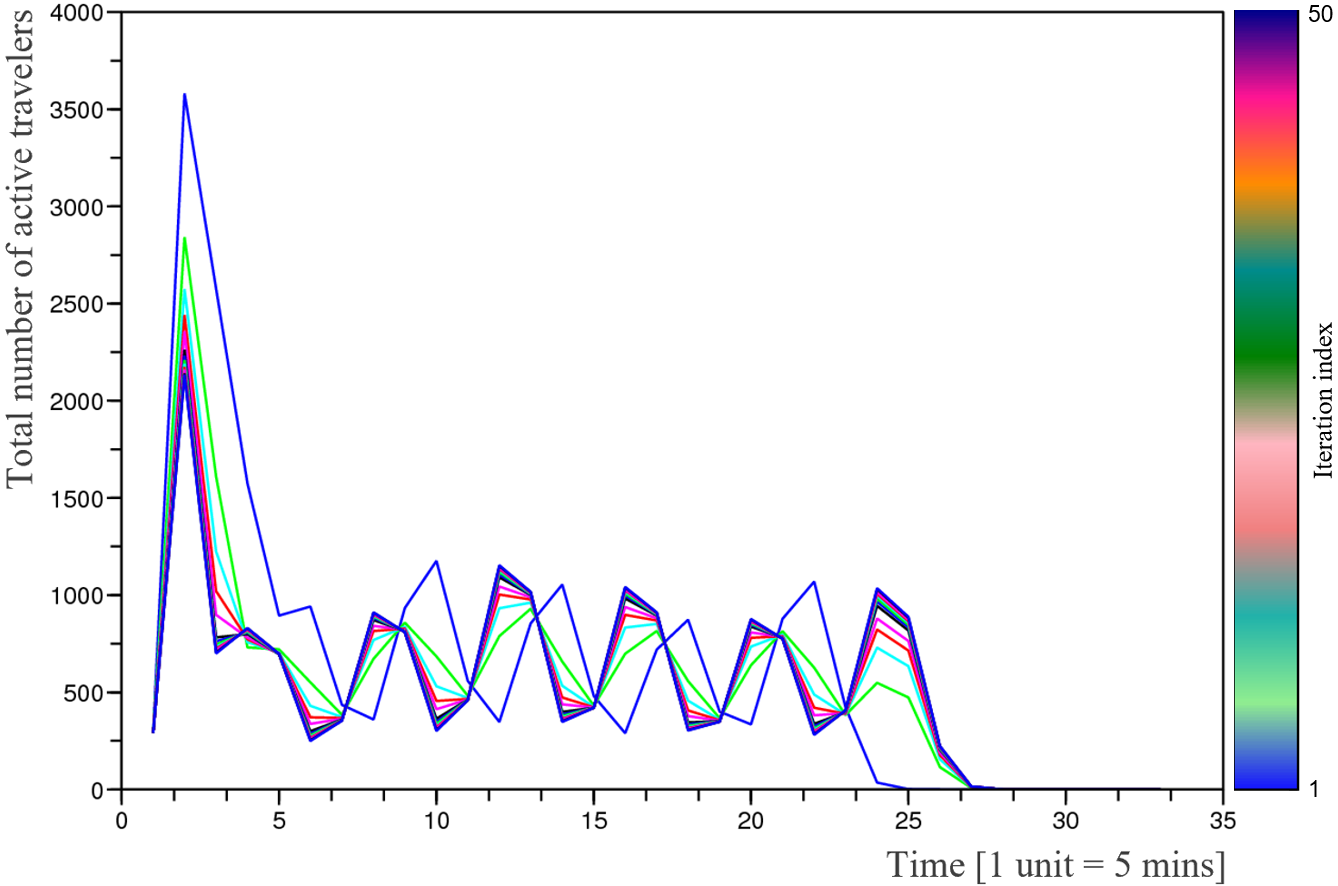}  
\caption{Convergence of the system measures during the optimization. On the left: Speed as a function of time and iteration. On the right: Total number of travellers in the system as a function of time and iteration.} 
\label{fig:itr1}
\end{center}
\end{figure}

Moreover, we investigate the evolution of the network criteria during the optimization. Figure \ref{fig:itr1} illustrates the results. Similar to the Paris test case, the proposed methodology follows the SO principles and converges to an SO solution smoothly. The convergence pattern illustrates that the solution method is computationally efficient as it converges with only ten iterations.

The oscillation observed in Figure \ref{fig:itr1} can be attributed to several factors: (i) The results are derived from the demand dataset, as shown in Figure 4, which itself exhibits a pronounced oscillatory feature. (ii) This oscillation is partially due to our specific definition of the desired arrival time. This definition leads to fluctuations in both average speed and vehicle accumulation, causing trips to cluster around the desired arrival times. Such clustering naturally forms V-shaped patterns around each desired arrival time value. In our previous study \cite{ameli2022departure}, which followed the approach of \cite{lamotte2018morning}, we distributed desired arrival times more evenly over time. This approach resulted in a smoother evolution of average speed and accumulation. (iii) Additionally, the penalty values in our model encourage trips to align as closely as possible with the desired arrival times, further contributing to this trend.

\begin{figure}[!h]
\begin{center}
\includegraphics[width=14cm]{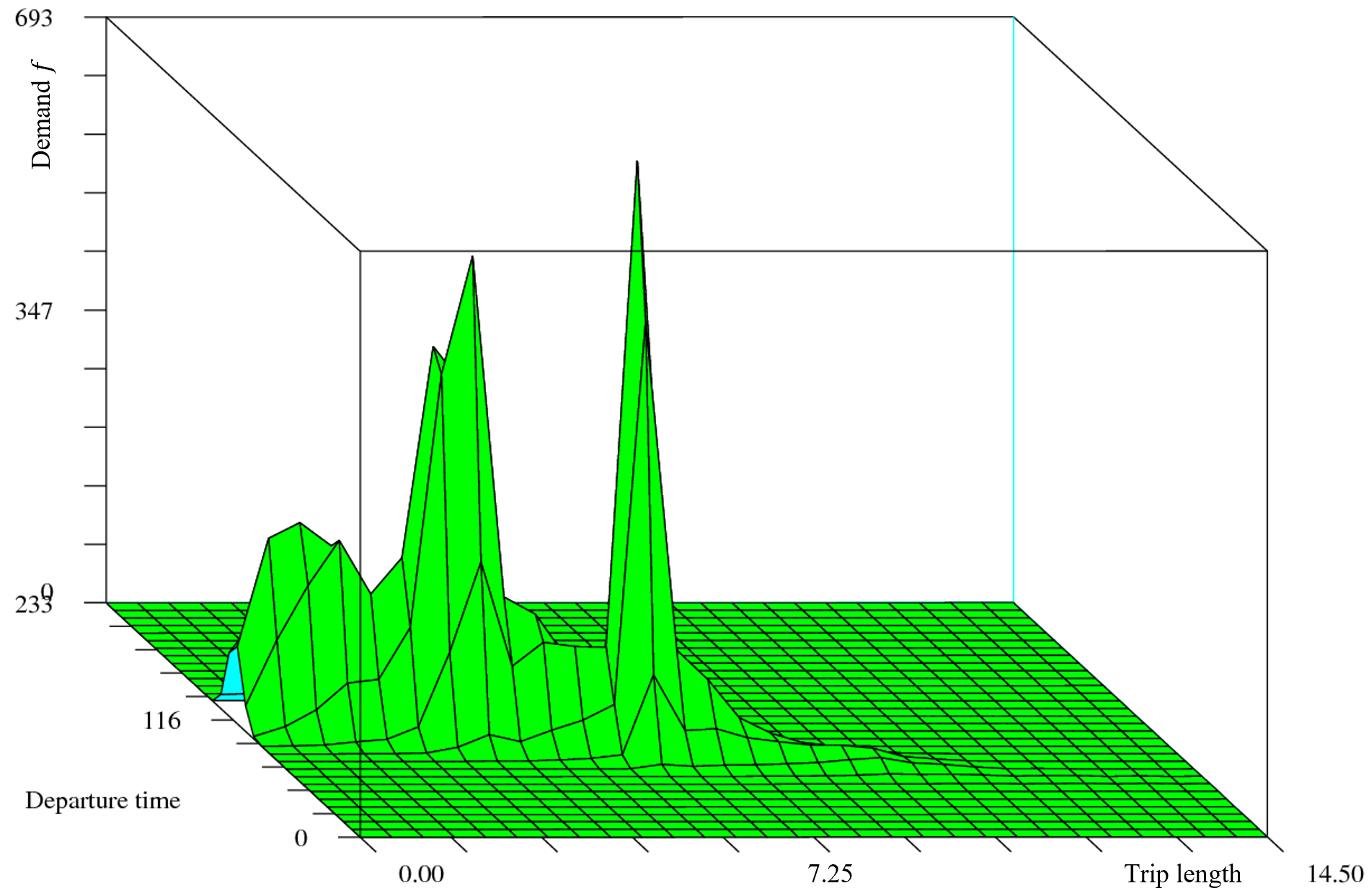}  
\caption{The density of departure times $f$ for a single arrival time (9 am).} 
\label{fig:Benchmark}
\end{center}
\end{figure}

The SO distribution of the departure time for the trips with the desired arrival time of 9:00 am is shown in Figure \ref{fig:Benchmark}. The results show that the solution for the SO does not follow any sorting pattern, e.g., FIFO and LIFO. In order to investigate further the solution characteristic compared to UE and SUE, we use a simulation-based framework for the large-scale full network of Lyon North with trip-based dynamic implementation.  

\subsection{Trip-based simulation for the SO problem}

The proposed methodology is extended for trip-based settings to represent commuters as a sort of agent with multiple attributes and decision variables \citep{zargayouna2008multi}. This subsection is structured to present the results of trip-based simulation on the Lyon North network and the benchmark analysis of three established network principles in dynamic traffic assignment: User Equilibrium (UE), Stochastic User Equilibrium (SUE), and System Optimum (SO). Since the proposed method can compute SO, we aim to evaluate these principles' solutions via simulation on a real test case. The trip-based simulator is designed using Particle discretization, presented in Section \ref{Sec: particle}.

\subsubsection{Validation of trip-based simulator}

Our simulator only keeps the original trip lengths as the generalized bathtub model does not account for the local traffic dynamics. In the test case of Lyon North, some trips have origins or destinations outside the covered area (51,215 trips). It means their trip starts or ends not inside the region. For the next simulation, we will not consider these trips in the departure time optimization. Therefore, 11,235 trips are fully interior. For those, the original departure time is disregarded and a desired arrival time is assigned. We divide them into seven classes with different desired arrival times. Table \ref{tab:Demand}).

\begin{figure}[ht]
\centering
\includegraphics[width=8cm, height=5.2cm]{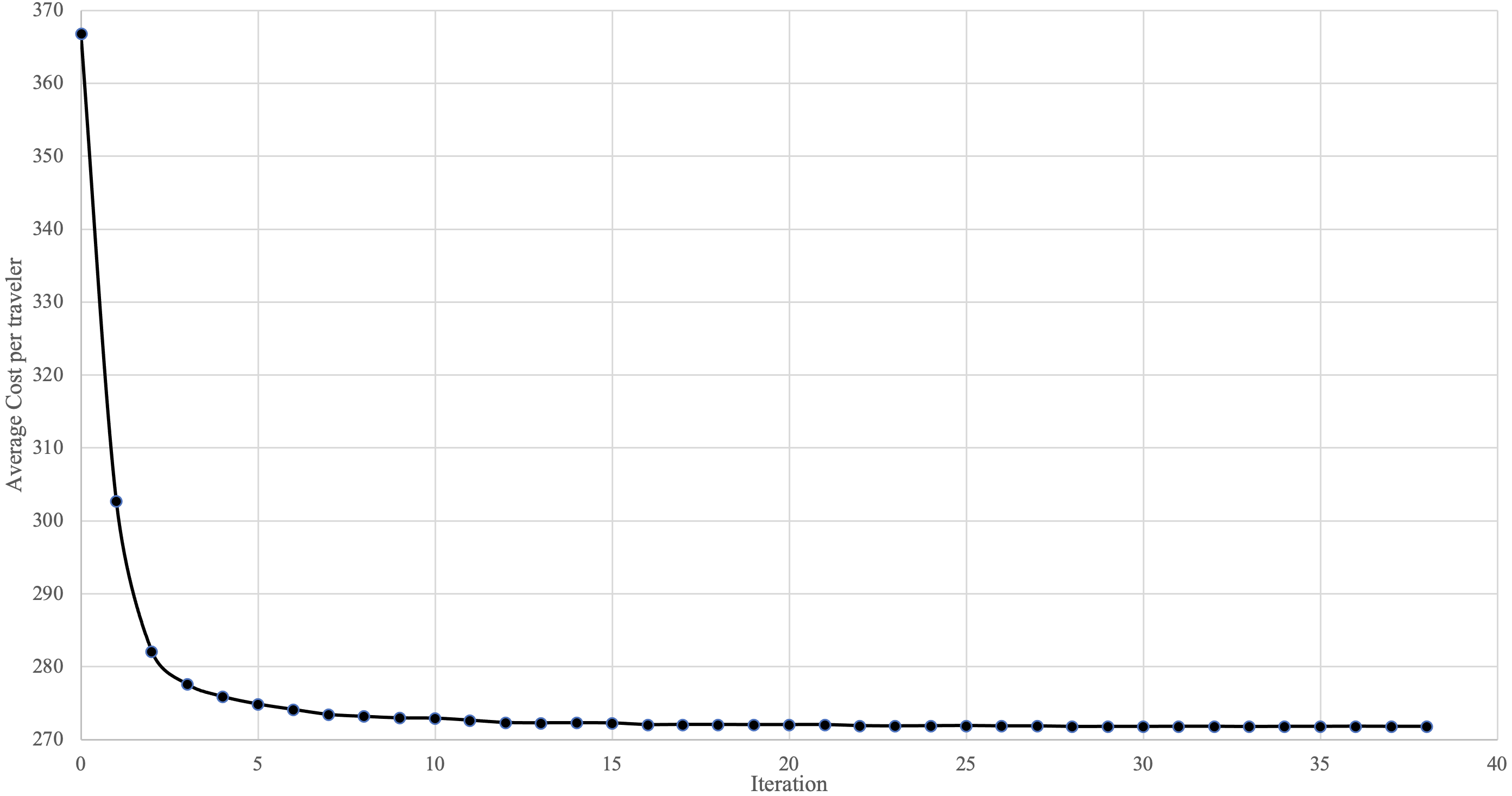} \includegraphics[width=8cm]{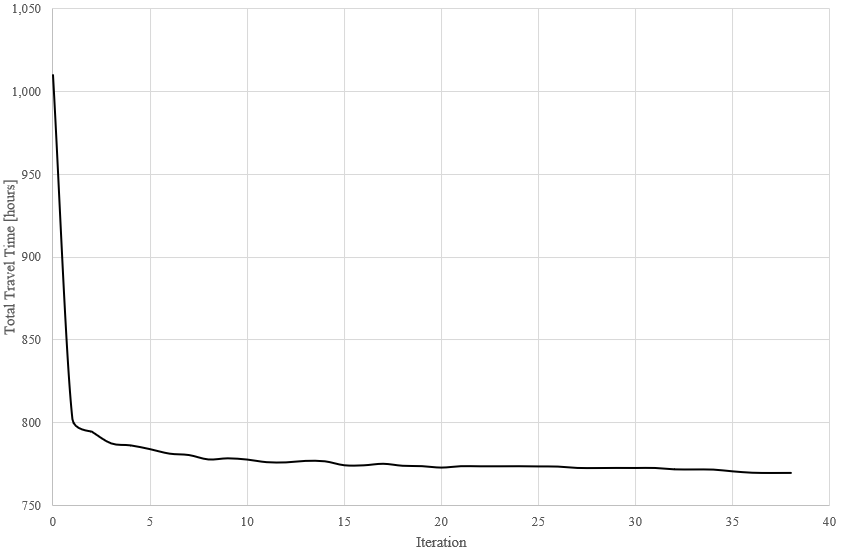}\hspace{5pt}
\caption{Evolution of Average Cost and Total Travel Time for targeted trips in the optimization process. The average cost is calculated by dividing the objective function by the total number of targeted users.}
\label{fig:Convergence_SO}
\end{figure}

Figure \ref{fig:Convergence_SO} presents the convergence of the solution method. We use average cost per traveler and total travel time as the convergence indicators. The average cost is calculated by dividing the total cost (the objective function) by the total number of targeted users (11,235). The algorithm converges smoothly after the drastic drop at the beginning because the initial solution starts the process. Few vibrations can be observed in the convergence pattern of the total travel time, which can be justified by the discrete nature of this configuration.

\begin{figure}[ht]
\centering
\includegraphics[width=10cm]{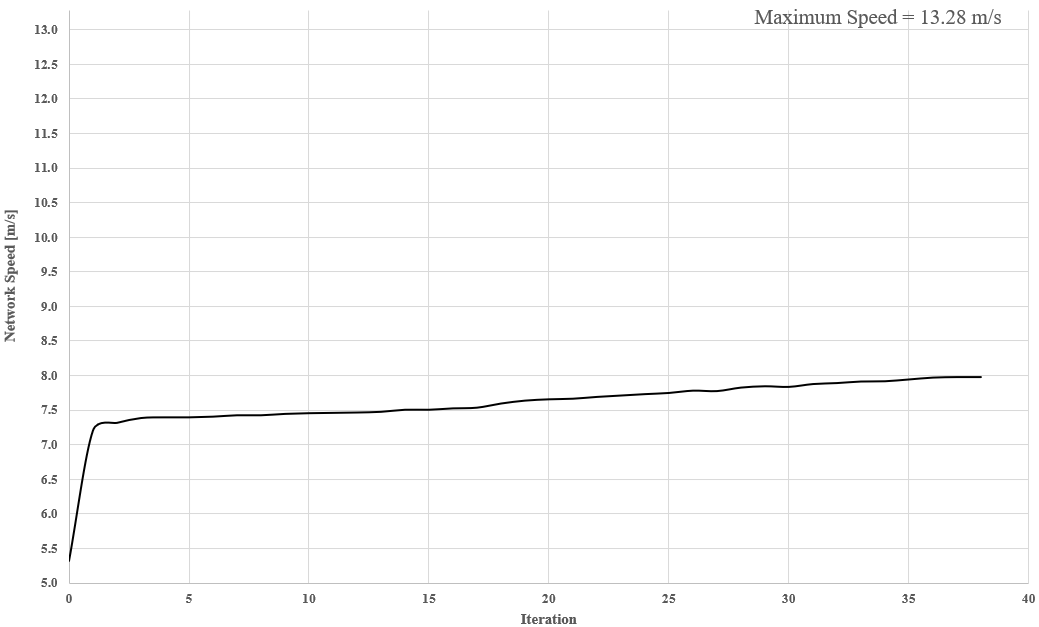}\hspace{5pt}
\caption{Evolution of Network Minimum Speed for targeted trips in the optimization process.}
\label{fig:Convergence_SO_speed}
\end{figure}

In addition, we present the evolution of the minimum network speed throughout each simulation during the optimization process to track the convergence of the algorithm and stability of the final solution. As depicted in Figure \ref{fig:Convergence_SO_speed}, the algorithm converges with few iterations, resulting in a final solution that demonstrates consistent characteristics. This implies that the gradient can no longer significantly improve the solution by minimizing the objective, indicating the stability of the final solution provided by the algorithm. Note that the minimum network speed denotes the minimum speed resulting from the speed function during the whole simulation (4 hours and 10 minutes) at every iteration. The network free-flow speed is equal to 47.8 km/h ($v_{max} = 13.28 m/s$), which is a standard value for a city-scale network. 

The simulation results show the consistency of our discretization method to capture the network dynamics in large-scale urban areas.

\subsubsection{Equilibria benchmark}

Since the simulation results for SO calculation is stable, we carried out the calculations of UE and SUE in order to compare the solutions of all three principles based on the network performance and trips indicators. 

The optimization process for UE and SUE is different. We calculate UE based on Mean Field Games framework that we previously developed in \cite{ameli2022departure}. While the SUE solution is calculated by the $f$ method presented in \cite{ameli2023morning, lebacque2022stochastic}. For all equilibria, iterative algorithms are applied. Each algorithm is started with the same initial solution where the targeted travelers with a higher trip length in all classes start their trip sooner than others based on the network free-flow speed ($v_{max} = 13.28 m/s$). The algorithms converge after $56$ iterations for UE and $21$ iterations for SUE to an equilibrium approximation.

\begin{figure}[ht]
\centering
\includegraphics[width=0.95\textwidth]{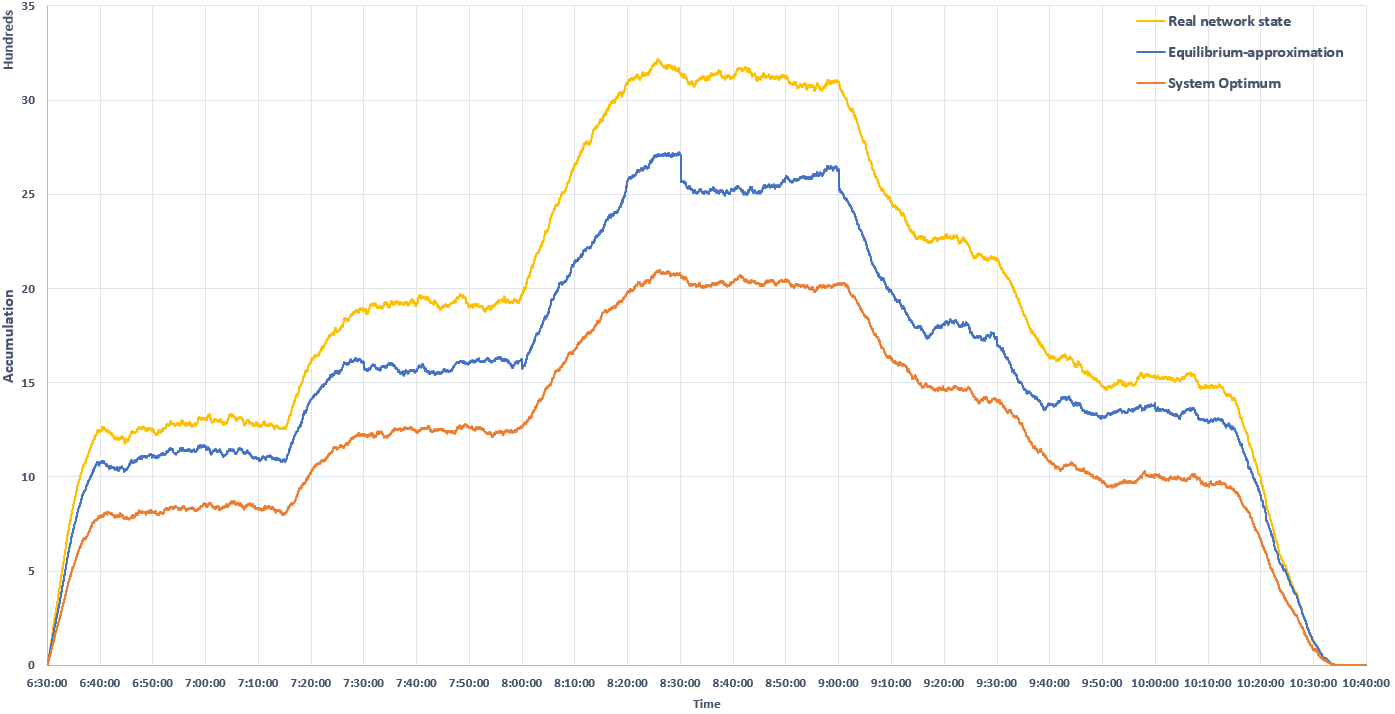}\hspace{5pt}
\caption{Results of the network’s performance overtime ($\Delta t = 1s$) in the different states: Accumulation of the real state of the network versus user equilibrium approximation and system optimum for all trips}
\label{fig:SO_Accumulation}
\end{figure}

Figure \ref{fig:SO_Accumulation} presents the equilibrium accumulation for the full demand, including targeted trips and background traffic, at each time step ($\Delta t = 1 \: sec$). It means the accumulation evolution in this figure is drawn for all trips, including exteriors that impact the network dynamics. The figure also includes the cumulative time series corresponding to the initial demand patterns with all given departure times. This curve sits above the UE and SO curves. Hence, the solution offered by the UE potentially enhances the cumulative travel time incurred by all users in the system in real network scenarios. The space between the cumulative time series determines this improvement. The SO accumulation is located below the UE, which could be expected because the SO minimizes the total cost, not necessarily the total travel time. Therefore, minimization of total travel time by SO is not necessarily expected. The final result depends on the desired arrival time and the cost of early and late penalties. 

To assess the deviation from the desired arrival times and corresponding early and late penalties, we grouped users within each desired arrival time category according to their trip lengths, using a 50-meter interval for segmentation. This approach resulted in approximately 58-87 clusters for each category. For each cluster, we calculated the travel cost difference by comparing the cost incurred by each user against the minimum cost within that cluster. Then, we normalized this difference against the minimum cost. This calculation aims to measure cost differences at the SO solution for users with similar arrival times and trip lengths (i.e., within the same cluster). For the UE solution, this measure, ideally, is approaching zero.

\begin{figure}[ht]
\centering
\includegraphics[width=0.75\textwidth]{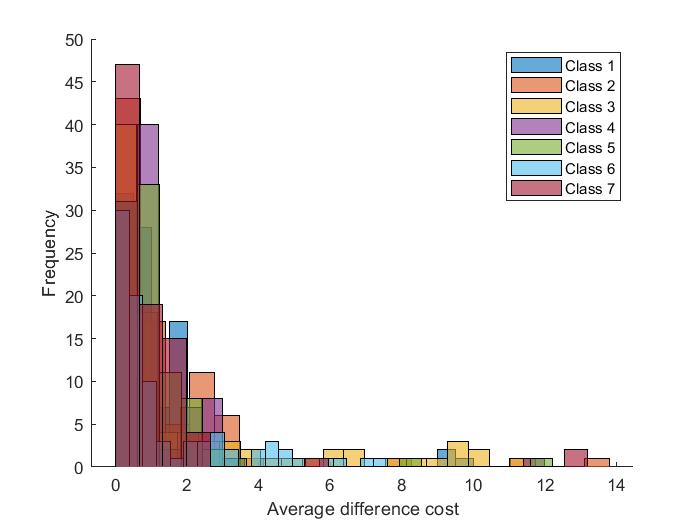}\hspace{5pt}
\caption{Optimization results regarding the different classes of trips. Note that the interior users are in the optimization process.}
\label{fig:Quality}
\end{figure}

Figure \ref{fig:Quality} presents the results of this measure for the SO solution. The findings reveal that the solution derived through our framework approximates a solution with error margins at the user level, which can be expected as the proposed framework calculates the SO solution. Notably, over 65\% of users achieved the optimal cost. A deviation of 7\%, was observed mainly among users in the late peak hours (classes 5-7).

Table \ref{tab:NetI} presents an analysis of network performance indicators for the three equilibria: SO, UE, and SUE. It provides a comparison of performance across these states based on several key indicators. The first indicator analyzed is Total Cost. Given that this indicator essentially represents the objective function of SO, it is expected to be the lowest for this state. Following the pattern observed for total travel time, the UE state's average travel cost is less than that of SUE.
The second indicator is Total Travel Time. As anticipated, the SO state shows the least total travel time compared to UE and SUE, with the highest travel time being attributed to the SUE. This is predictable considering SUE reflects the biases in commuter decision-making arising from their imperfect knowledge to calculate perceived costs. The third indicator analyzed is Average Cost, calculated by dividing the total cost by the total number of targeted users (11,235). 

To illustrate the variability in travel costs, the table includes the standard deviation (STD) of this indicator. Notably, the SO state exhibits the highest variance in travel cost, aligning with the idea that this principle favors overall system performance over individual user gains. Interestingly, the STD for SUE is less than SO but greater than UE, though closer to the UE value. This can be explained by the parameters of the logit function used in the SUE model.

The final indicator, Average Delay, presents significant differences between the equilibria. It is determined by the absolute difference between the actual arrival time and the desired arrival time without any penalties. To calculate the average delay, we divide the total delay by the total number of targeted users. As might be expected, the SO solution has the highest average delay value. The UE state's average delay, however, is notably lower than both SO and SUE. This discrepancy arises from the fact that the average delay can mirror to the objective function of UE and is correlated with total travel cost. This variation also exemplifies the 'price of anarchy' observed between the SO and the UE/SUE solutions.

\begin{table}[]
\centering
\caption{Network performance indicators for the equilibria}
\begin{tabular}{c|c|c|c|c|c}
\multirow{2}{*}{\begin{tabular}[c]{@{}c@{}}Equilibrium \\  Principle\end{tabular}} & \multirow{2}{*}{\begin{tabular}[c]{@{}c@{}}Total travel \\ cost \end{tabular}} & \multirow{2}{*}{\begin{tabular}[c]{@{}c@{}}Total travel time\\ (hours)\end{tabular}} & \multicolumn{2}{c|}{Average cost} & \multirow{2}{*}{\begin{tabular}[c]{@{}c@{}}Average delay\\ (min)\end{tabular}} \\ \cline{4-5}
                                         &                                                                                      &  & Mean      & Standard deviation   &                                                                                \\\hline
UE                &           3672946,20            & 930.02                                                                               & 326.92    & 23.95                & 0.39                                                                           \\
SUE               &         3725350,77           & 1247.61                                                                              & 331.58    & 24.95                & 1.48                                                                           \\
SO                &           3053335,95            & 769.9                                                                                & 271.77    & 38.55                & 7.09                                                                          
\end{tabular}\label{tab:NetI}
\end{table}

\begin{table}[ht]
\centering
\caption{Demand profile and the results for multi-class users of Lyon North}
\resizebox{\textwidth}{!}{%
\begin{tabular}{c|c|c|c|c|c|c|c|c|c|c|c}
\multirow{2}{*}{Class} & \multirow{2}{*}{Share} & \multirow{2}{*}{Number} & \multirow{2}{*}{\begin{tabular}[c]{@{}c@{}}Mean trip\\ length (km)\end{tabular}} & \multirow{2}{*}{\begin{tabular}[c]{@{}c@{}}Arrival\\ time\end{tabular}} & \multirow{2}{*}{\begin{tabular}[c]{@{}c@{}}Desired\\ arrival time\end{tabular}} & \multicolumn{3}{c|}{Average cost} & \multicolumn{3}{c}{Average delay (min)} \\ \cline{7-12} 
                                      &                                     &                                      &                                        &                                        &                                                   & UE        & SUE       & SO       & UE          & SUE          & SO       \\ \hline
Class 1 & 13.73\% & 1,543 & 2.53 & 6:30-7:15 & 7:00 & 294.7738 & 393.0646 & 259.1457 & 0.4994 & 1.9014 & 5.9841 \\
Class 2 & 13.84\% & 1,555 & 2.58 & 7:15-7:45 & 7:30 & 309.7068 & 323.4421 & 278.7025 & 0.4206 & 1.5863 & 7.5002 \\
Class 3 & 15.42\% & 1,732 & 2.55 & 7:45-8:15 & 8:00 & 313.0173 & 318.8728 & 309.1616 & 0.3355 & 1.2551 & 8.9018 \\
Class 4 & 18.30\% & 2,056 & 2.65 & 8:15-8:45 & 8:30 & 381.8035 & 322.2873 & 286.9059 & 0.4130 & 1.5484 & 7.8546 \\
Class 5 & 15.05\% & 1,691 & 2.63 & 8:45-9:15 & 9:00 & 358.3773 & 303.3978 & 259.1457 & 0.3228 & 1.2337 & 7.3683 \\
Class 6 & 11.82\% & 1,328 & 2.70 & 9:15-9:45 & 9:30 & 313.9684 & 299.1294 & 241.1212 & 0.2997 & 1.1226 & 5.7597 \\
Class 7 & 11.84\% & 1,330 & 2.63 & 9:45-10:30 & 10:00 & 290.5534 & 368.9470 & 242.7808   & 0.4402 & 1.7113 & 5.3444 \\
\end{tabular}\label{tab:Demand}
}
\end{table}

Table \ref{tab:Demand} provides an overview of each user class considered for the optimization process. It delineates the distinctive characteristics of each class, such as their proportion of the total demand, average trip length, and desired arrival times, indicating individual scheduling habits.

Focusing on commuters, Table \ref{tab:Demand} outlines key indicators, including the average cost and average delay associated with every class. Crucially, this data is presented across three equilibria for a comparative analysis.

Interestingly, the average cost for all user classes in the SO equilibrium is consistently lower than that in both the UE and SUE equilibria, aligning with the aggregate results presented in Table \ref{tab:NetI}. However, the relationship between the SUE and UE is not direct. For instance, the average cost for SUE in the first three user classes exceeds that of UE, yet it subsequently falls below the UE value before rising again. This trend indicates that the UE solution presents a superior average cost compared to SUE between the peak-hour window of 7:30 am to 9:30 am.

In terms of average delay, the relationship between the three equilibria mirrors that of the collective values. The results corroborate the finding that at high congestion levels, the SO and UE states differ substantially in terms of travel delays, further justifying the magnitude of discrepancy between the network performance of SO and UE. 

\section{Conclusion and future works} \label{sec:conclusion}

This study presents a novel formulation for the departure time system optimum problem based on the generalized bathtub model, providing a generic approach to capture the complex dynamics of urban traffic congestion. By incorporating a continuous formulation that can accommodate any distribution for trip length and desired arrival time, the proposed framework offers a more realistic representation of the heterogeneous characteristics of trips in an urban setting. The method can be extended to the case when a downstream supply constraint is present.

The application of the proposed methodology to the morning commute problem of the Lyon North network demonstrates its effectiveness in solving the system or social optimum (SO) problem with multiple desired arrival times and heterogeneous trip lengths for a large number of trips. The existence of the SO solution is proven. However, the conclusion regarding the uniqueness of the solution is not trivial. Further investigations in this direction should assess the influence of the structure of the travel demand density $m(t_a,x)$ on the presence of local optima. Additionally, an analytical process is introduced to calculate the marginal travel time for solving the SO problem, enhancing computational efficiency.

Furthermore, a benchmark analysis comparing the solution of User Equilibrium (UE), Stochastic User Equilibrium (SUE), and SO shows that the proposed methodology outperforms UE and SUE solutions in terms of network performance indicators, specifically total travel cost (time), which was expected based on the definition of these principles.

The authors could outline several future research directions, including conducting an analytical test case to further investigate the features of the continuous model. Additionally, comparing the results of different discretization approaches and benchmarking the model with other equilibrium models for macroscopic and microscopic models are ongoing efforts. These endeavors aim to strengthen the understanding and applicability of the proposed framework in addressing urban traffic congestion.



\section*{Authors contribution statement}

\textbf{MA}: Conceptualization;  Formal analysis; Investigation; Methodology; Software; Project administration; Visualization; Writing - original draft; Writing - review \& editing. \textbf{JPL}: Conceptualization;  Formal analysis; Investigation; Methodology; Software; Project administration; Visualization; Writing - original draft; Writing - review \& editing. \textbf{NA}: Conceptualization;  Formal analysis; Investigation; Methodology; Software; Project administration; Visualization; Writing - original draft; Writing - review \& editing. \textbf{LL}: Conceptualization;  Formal analysis; Methodology; Project administration; Visualization; Writing - original draft; Writing - review \& editing.

\section{Acknowledgements}

This research received no specific grant from funding agencies in the public, commercial, or not-for-profit sectors.

\section*{Conflicts of interest}
None.

\end{sloppypar} 

\typeout{}
\bibliography{references}

\newpage

\appendix

\setcounter{table}{0}
\setcounter{figure}{0}

\end{document}